\documentclass[10pt]{article}
\input epsf.sty

\epsfverbosetrue %messages will show height and width

\newcommand{\xincludeps}[3]{\begin{figure}
                                 \centerline{%
                                 \epsfxsize=#2in
                                 \epsfbox{#1}
                                 }
                            \caption{ {\sf #3}
\label{fig:#1}}
                            \end{figure}
                           }

\usepackage{latexsym}
\usepackage{amsfonts}
\usepackage{amscd}

\newcommand{\R}{\mathbb{R}}

\newcommand{\Z}{\mathbb{Z}}
\newcommand{\C}{\mathbb{C}}
\newcommand{\IP}{\mathbb{P}}
\newcommand{\IH}{\mathbb{H}}

\def\cbar{\overline{\C}}
\newcommand{\uu}{\ensuremath{{\bf{u}}}}
\newcommand{\ww}{\ensuremath{{\bf{w}}}}

\newcommand{\rs}{\mbox{$\widehat{\C}$}}

\def\DDD{{\mathcal D}}

\def\FFF{{\mathcal F}}

\def\HHH{{\mathcal H}}

\def\JJJ{{\mathcal J}}
\def\LLL{{\mathcal L}}
\def\MMM{{\mathcal M}}

\def\MMM{{\mathcal M}}
\def\WWW{{\mathcal W}}

\def\NNN{{\mathcal N}}

\def\UUU{{\mathcal U}}

\newtheorem{thm}{Theorem}[section]
\newtheorem{defn}{Definition}[section]

\newtheorem{prop}{Proposition}[section]

\newtheorem{lemma}{Lemma}[section]
\newtheorem{cor}{Corollary}[section]

\newcommand{\qed}{\nopagebreak
\begin{flushright}%end-of-proof
         \rule{2mm}{2.5mm} \end{flushright}}

\newcommand{\bdry}{\partial}                     %boundary
\newcommand{\id}{\mbox{\rm id}}                  %identity
\newcommand{\cl}{\overline}                      %closure
\newcommand{\Imag}{\mbox{\rm Im}}                    %Im, in better font
\newcommand{\Real}{\mbox{\rm Re}}                    %Re, in better font

\newcommand{\Aut}{\mbox{\rm Aut}}
\newcommand{\Mod}{\mbox{\rm Mod}}    %Mod
\newcommand{\Teich}{\mbox{\rm Teich}}%Teich
\newcommand{\intersect}{\cap}
\newcommand{\union}{\cup}                        %union
\newcommand{\mtwo}[4]                            %2x2 matrices--
{\mbox{$\left(\begin{array}{cc}                  %takes four arguments #1
#1& #2 \\ #3 & #4
\end{array}
\right)$}}

\newcommand{\dettwo}[4]                          %2x2 matrices--
{\mbox{$\left|\begin{array}{cc}                  %takes four arguments #1
& #2 \\ #3 & #4
\end{array}
\right|$}}

\newcommand{\pf}{\noindent {\bf Proof: }}

\newcommand{\be}{\begin{enumerate}}
\newcommand{\eb}{\end{enumerate}}
\newcommand{\bi}{\begin{itemize}}
\newcommand{\ib}{\end{itemize}}
\newcommand{\bl}{\begin{list}}
\newcommand{\lb}{\end{list}}

\newcommand{\gap}{\vspace{5pt}}             %make a space of a blank line

\newcommand{\Rat}{\mbox{\rm Rat}}

\newcommand{\Fix}{\mbox{\rm Fix}}

\newcommand{\genby}[1]{\mbox{$\langle #1 \rangle$}}

\newcommand{\Crit}{\mbox{\rm Crit}}

\newcommand{\Poly}{\mbox{\rm Poly}}  %Poly
\newcommand{\GRat}{\mbox{\rm GRat}_d}
\newcommand{\GRatmn}{\mbox{\rm GRat}^{\times, *}_d}
\newcommand{\GRatm}{\mbox{\rm GRat}^{\times}_d}

\begin{document}

\title{Spinning deformations of rational maps}

\author{ Kevin M. Pilgrim\footnote{Research partially supported by NSF
grant No. DMS 9996070.}  \  and Tan Lei}

\maketitle

\tableofcontents 

\newpage

\section{Introduction}
\label{secn:results}

We analyze a real one-parameter
family of quasiconformal  deformations of a hyperbolic rational map known
as {\em spinning}.  We show that under fairly general hypotheses, the limit
of spinning either exists and is unique, or else converges to infinity in
the moduli space of rational maps of a fixed degree.  When the limit
exists, it has either a parabolic fixed point, or a pre-periodic
critical point in the Julia set, depending on the combinatorics of the
data defining the deformation. The proofs are  soft and rely on two
ingredients:  the construction of a Riemann surface 
containing the closure  of the family, and an analysis of the
geometric limits of some simple dynamical systems.  
 An interpretation in terms of Teichm\"uller theory is presented as well.

\subsection{Definition of spinning} Although we shall only treat the
simplest cases in this work, we define spinning in a general context.
\gap

\noindent{\bf Notation.}  Let $f: \IP^1 \to \IP^1$ be a rational map of
degree $d \geq 2$.  The {\em grand orbit} of a point $z \in \IP^1$ is the
set of $w$ such that
$f^{\circ i}(z) = f^{\circ j}(w)$ for some $i,j
\geq 0$. Let
$\hat{J}$ denote the closure of the grand orbits of all periodic points and
all critical points of $f$.  The set
$\hat{J}$ contains the Julia set.  Its complement
$\hat{\Omega}$ is therefore contained in the Fatou set and is the disjoint
union of open subsets $\Omega^{dis}
\sqcup
\Omega^{fol}$ consisting of points whose grand orbits are discrete and
indiscrete sets, respectively.   For generic hyperbolic rational maps,
$\Omega^{fol}$ is empty.
\gap

\noindent{\bf Quotient surfaces.}  By the classification of stable regions,
the set $\Omega^{dis}$ consists of points lying in the basin of attraction
of parabolic and (non-super)-attracting cycles.  The restriction of $f$
to $\Omega^{dis}$ is a holomorphic self-covering, so the quotient
$\Omega^{dis}/f$ of $\Omega^{dis}$ by the grand orbit equivalence relation
is a one-dimensional, possibly disconnected complex manifold and will be
called the {\em quotient surface} associated to $f$. The quotient surface
is a disjoint, possibly empty union of at-least-once-punctured tori (one
for each attracting cycle) and at-least-once-punctured copies of $\C^*$
(one for each parabolic cycle), where in each case the number of
punctures is the number of grand orbits of critical points in the
corresponding basin.  Note that by a theorem of Fatou, the {\em
immediate} basin of an attracting or parabolic cycle contains at least
one critical point, and yields therefore at least one puncture in each
case.

Let $\Psi:
\Omega^{dis} \to
\Omega^{dis}/f$ be the canonical projection.

\gap

\noindent{\bf Input data for spinning.} 
Let $S$ be a component of $\Omega^{dis}/f$ with
 at least two punctures. Let $\Psi(c)$ be one of these
punctures, where $c$ is a critical point of $f$ (abusing notation, we write
$c=\Psi(c)$ if no confusion can arise). The point $c$ we call 
the {\em spun critical point}.  Let
$[\gamma] \in \pi_1(S^\sharp, c)$ be a homotopy class of oriented simple
closed curve on the topological surface
$S^\sharp = S\union\{c\}$ passing through $c$.  

To describe spinning as a continuous, as opposed to a discrete, process
we make the following non-canonical choices. 

Let $A$ be an annulus in $S^\sharp$
with (oriented) core curve
$\gamma$ and real-analytic boundary (so that $A$ does not contain the other punctures).
Choose a universal cover $\widetilde{A} = \{x+iy \in \C | -2l < y < 2l
\}$  with
$p_A:\widetilde{A} \to A$  the  covering map, such that $p_A(\R)=\gamma$, $p_A$
 is orientation-preserving on $\R$, $p_A(0)=c$, and that
 the map on $\widetilde{A}$ induced by $\gamma$ is translation by one
(to the right).   This implies that the modulus of $A$ is
$4l$. With these conventions, the covering space and map
$p_A: \widetilde{A} \to A$ are unique.   
\gap
 
\noindent{\bf Definition of spinning.}  Let $f_0 = f$.    
Define
\[ \widetilde{h}: \widetilde{A} \times \R \to \widetilde{A}\] by linearly
interpolating translation to the right by $t$ on the horizontal strip $|y|
\leq l$ and the identity map on the boundary of
$\widetilde{A}$.  More precisely, set
\[ \widetilde{h}(x+iy, t) =
\left\{
\begin{array}{lcl} x+t+iy               & \mbox{ if } &  0
\leq |y| \leq l,\\  x+t(2-|y|/l) +iy     & \mbox{ if } &  l
\leq |y| \leq 2l.
\end{array}
\right.
\] We set $\widetilde{h}^t =\widetilde{h}(\cdot, t):
\widetilde{A} \to
\widetilde{A}$.  The key features are:
\bi
\item For each $t$, the map $\widetilde{h}^t$ is quasi-conformal, 
$\widetilde{h}^t|_{\partial \widetilde{A}}=id$,
 and $\widetilde{h}^t$ commutes with the group of deck
transformations, thus giving a well-defined map $h: A \times \R \to A$.

\item $h$ extends to a continuous homomorphism from the reals under
addition to the group of qc self-homeomorphisms of the Riemann surface
$S^\sharp$, such that $h^t$ is the identity on the complement of $A$.
If e.g. $S$ is a torus $T$ with punctures, then $h$ extends to $T$ as
well.

\item For each $t$, the map $h^t: S^\sharp \to S^\sharp $ is conformal on
the complement of the region
$p_A( l< |y| < 2l )$, a union of two parallel subannuli of
$A$.  Also $h^t(A)=A$, so that the modulus is unchanged.  
\ib

\xincludeps{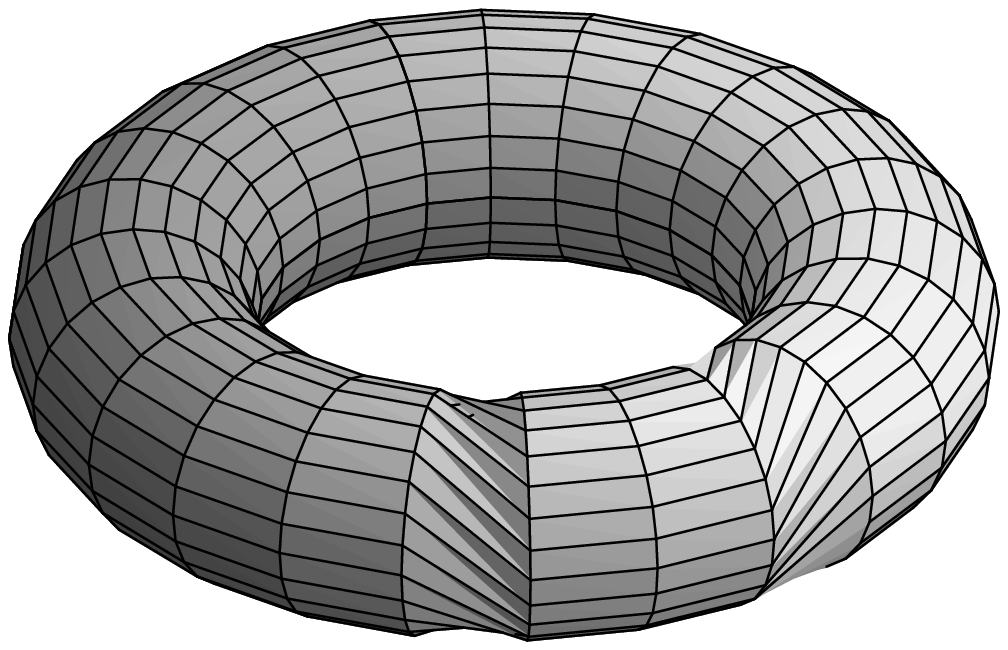}{3.0}{Spinning on the torus.}

Finally, let $\mu_t$ be the $f$-invariant Beltrami differential obtained by
pulling back the dilatation of
$h^t$ under the canonical projection
$\Psi: \Omega^{dis} \to \Omega^{dis}/f$.  By the Measurable Riemann
Mapping Theorem there is a quasiconformal homeomorphism
\[ H_t: \rs \to \rs, \] unique up to postcomposition with M\"obius
transformations, whose dilatation agrees with
$\mu_t$.  If a representative $H_t$ is chosen, then the map
\[ f_t = H_t \circ f \circ H_t ^{-1}\] preserves the standard conformal
structure and is therefore a rational map.  Different representatives for
$H_t$ yield M\"obius conjugate
$f_t$, and so we obtain a map
\[ \sigma = \sigma_{\gamma,A}: \R \to \Rat_d/\Aut(\IP^1)\] from the reals into
the space of conjugacy classes of rational maps of degree $d$, which we
call a {\em spinning path} of
$\gamma$.   The image
$\sigma([0, \infty))$ we call a {\em spinning ray} of
$\gamma$.  It is easily shown that $\sigma$ is real-analytic.  
\gap

\noindent{\bf Visibility.}  A curve $\gamma \subset S^\sharp$ can lift to
the sphere in a variety of combinatorially distinct ways.   Let $S,
\gamma$  be as in the setup for spinning, where
$S=\Psi(\hat{B})$ and $B$ is the immediate basin of an attracting 
periodic cycle.   Let
\[ \Gamma=\bigcup_{\delta \subset \Psi^{-1}(\bdry A)}
\cl{\delta} \] where the union is over all connected components of the
preimage of
$\bdry A$ under the projection $\Psi$.  Then $\Gamma$ divides $B$ into
various components; see Figures
\ref{fig:vis.eps}, \ref{fig:invis.eps} where $\Psi^{-1}(A)$ consists of
the lighter colored regions.  Observe that if 
$W_1, W_2, ...$ denote the components of $B-\Gamma$ whose closures
contain points of the attractor, then the $W_k$ are finite in number,
and the restriction $f^p: \union_k W_k \to \union_k W_k$ is proper.

\begin{defn}[Visible point]
\label{defn:visible}  A point $w \in B - \Gamma$  is {\em visible with
respect to $\gamma$} if the closure of the unique connected component of
$B-\Gamma$ containing $w$ contains a point of the attractor.   A point $w
\in \IP^1$ becomes {\em visible with respect to $\gamma$ after $r$ steps}
if $f^{\circ r}(w)$ is visible but
$f^{\circ i }(w)$ is  not visible for $0 \leq i <r$.
\end{defn}

For example, the critical points $c$ and $b$ in Figure
\ref{fig:vis.eps} are both visible, while in Figure
\ref{fig:invis.eps} critical point $b$ is visible but $c$ is visible
after one step.   Note that a visible critical point
$b$  necessarily has infinite
forward orbit, hence $\phi(b) \neq 0$ in a linearizing coordinate $\phi$.

\xincludeps{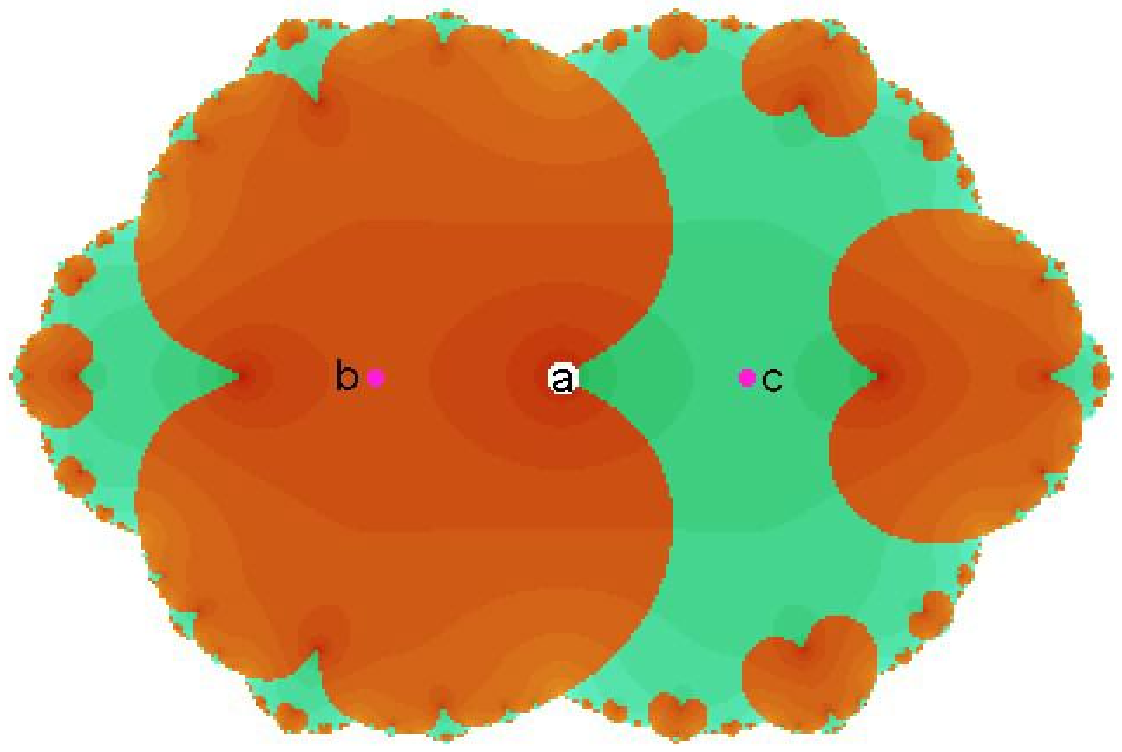}{3.0}{The critical point $c$ is visible.}

\xincludeps{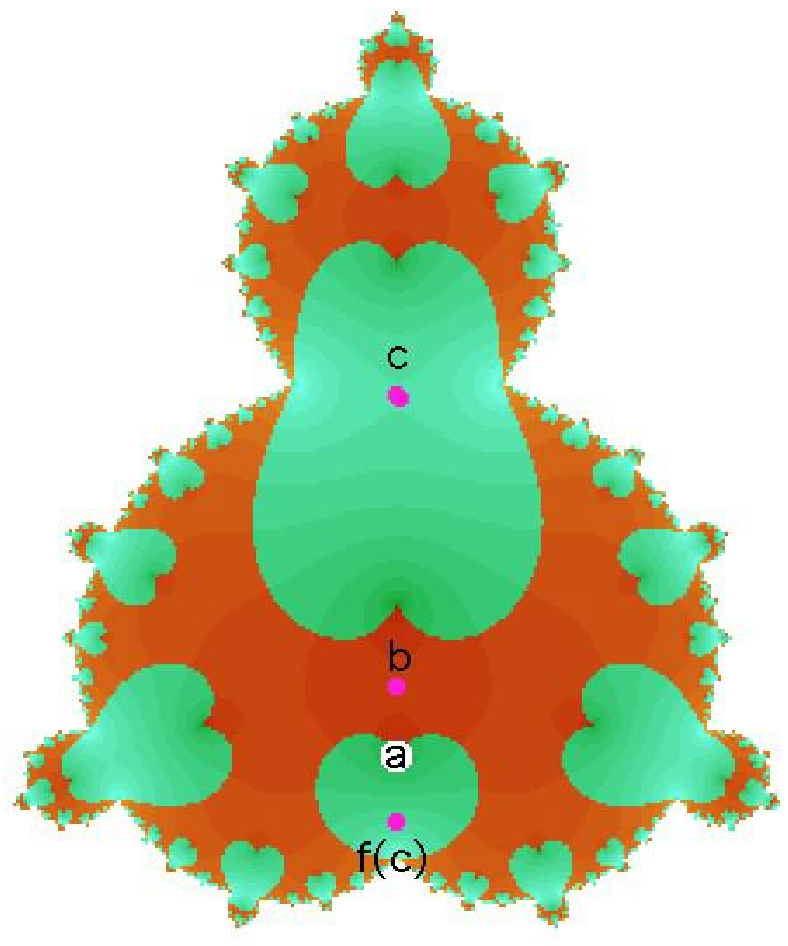}{3.0}{The critical point $c$ is not visible.}

\subsection{An example}

Consider, for complex $c \neq 0$, the family of cubic critically marked
polynomials
\[ f(c,z) = -\int_0 ^{z} (\zeta-c)(\zeta+\frac{1}{2c})d\zeta = \frac{1}{2}z
+ az^2 - \frac{1}{3}z^3\] where
$a=\frac{1}{2}(c-\frac{1}{2c})$.  This family consists of polynomials with
an attracting fixed point of multiplier
$1/2$ at the origin, normalized conveniently so the leading coefficient is
$-1/3$, and having two marked critical points $c$ and
$b=b(c)=-\frac{1}{2c}$.  Conjugation by $z \mapsto -z$ preserves this
family and replaces $c$ by $-c$.

Let us set
\[ c_0 = \frac{\sqrt{2}}{2} \approx .707.\] For $c=c_0$, the marked map
$f^\times_0(z) = f(c_0, z)$ is equal to
\[ f^\times_0(z) = -\frac{1}{3}z^3 + \frac{1}{2}z\] and in particular is
odd.  Hence both critical points $c_0$ and
$b_0=-c_0$ lie in the immediate basin of the origin, since by Fatou's
theorem the immediate basin must contain at least one critical point.

In the $c$-parameter plane, consider the locus
\[ \{ c \; | \; b \; \mbox{converges to the origin under iteration of
$f(c,z)$}\;\}.\]    

Let $X(f^\times_0)$ be the connected component containing
$f_0^\times$; see Figure
\ref{fig:parampair.eps}.   Then $X(f^\times_0)$ is the locus of maps
$f^\times_c$ such that the  critical point $b=b(c)$ lies in the immediate
basin of the origin.  Evidently $X(f^\times_0)$ is not closed in the
$c$-parameter plane.

Let $q,r$ be positive integers.  Let $X_{par}^{q, r}(f^\times_0)$ denote
the subset of $X(f^\times_0)$ for which the corresponding maps have a
parabolic cycle of period less than or equal to $q$ and of multiplier an
$r$th root of unity.  Let
$X_{mis}^{r, q}(f^\times_0)$ denote the set of Misiurewicz-type parameter
values, i.e. those for which the critical point $c$ lands on a repelling
cycle of period at most $q$  after at most
$r$ iterations. Note that 
$X_{par}^{q, r}(f^\times_0)$ and $X_{mis}^{q, r}(f^\times_0)$ are discrete
subspaces of $X(f^\times_0)$.

\xincludeps{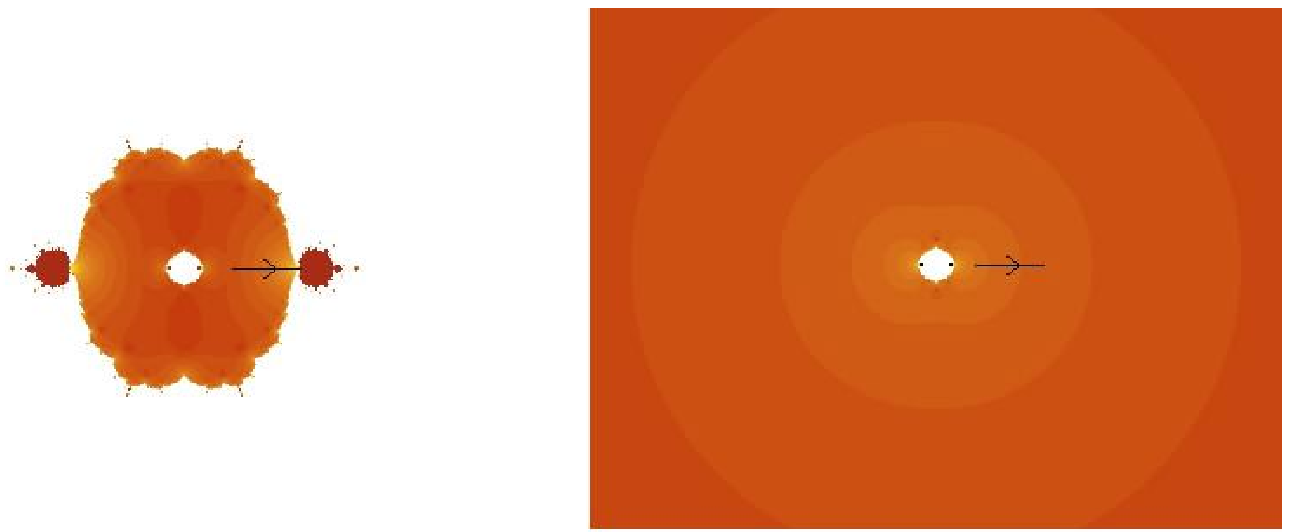}{4.0}{At left, the annular region shaded in
gray is the locus of $c$-parameters for which both critical points $c, b$
lie in the immediate basin of the attractor at the origin.  Superimposed
is the path $\sigma(t), t \geq 0$.   At right, the region shaded in gray
is the locus $X(f^\times_0)$ of $c$-parameters for which the critical
point $b$ lies in the immediate basin of the attractor at the origin. 
Both images are the windows $|\Real(c)| \leq 5.6, |\Imag(c)| \leq 4.2$. 
Note that the limit of $\sigma(t)$ lies in $X(f^\times_0)$.}

Let $\gamma$ be a simple closed curve satisfying the hypothesis of Theorem
\ref{thm:lands_or_diverges}.  For example, we may take $\gamma$ to  be the
image of the positive real axis under the projection $p:
\C^* \to T$, i.e. a real curve containing $\Psi(c)$.  Choose an annulus
$A$ surrounding $\gamma$, and consider spinning for the map $f_0$ about
the curve $\gamma$.   By examining the real graph of $f_0$ it is easy to
see that both $b$ and $c$ are  visible with respect to $\gamma$; cf
Figure
\ref{fig:vis.eps} in which $\Psi^{-1}(A)$ for
$f_0$ is shown as the lighter colored region.

A lift $\sigma^\times $ of the spinning ray to $X(f^\times_0)$ is
shown in Figure \ref{fig:parampair.eps}.  The evolution of the dynamics
is shown in Figure \ref{fig:evoln4.eps}.    The limiting map in the
above example is
\[ g^\times(z)=f(c_\infty, z)\] where
\[ c_\infty =\frac{\sqrt{2}+\sqrt{14}/2}{\sqrt{3}}
\approx 1.897.\]

\xincludeps{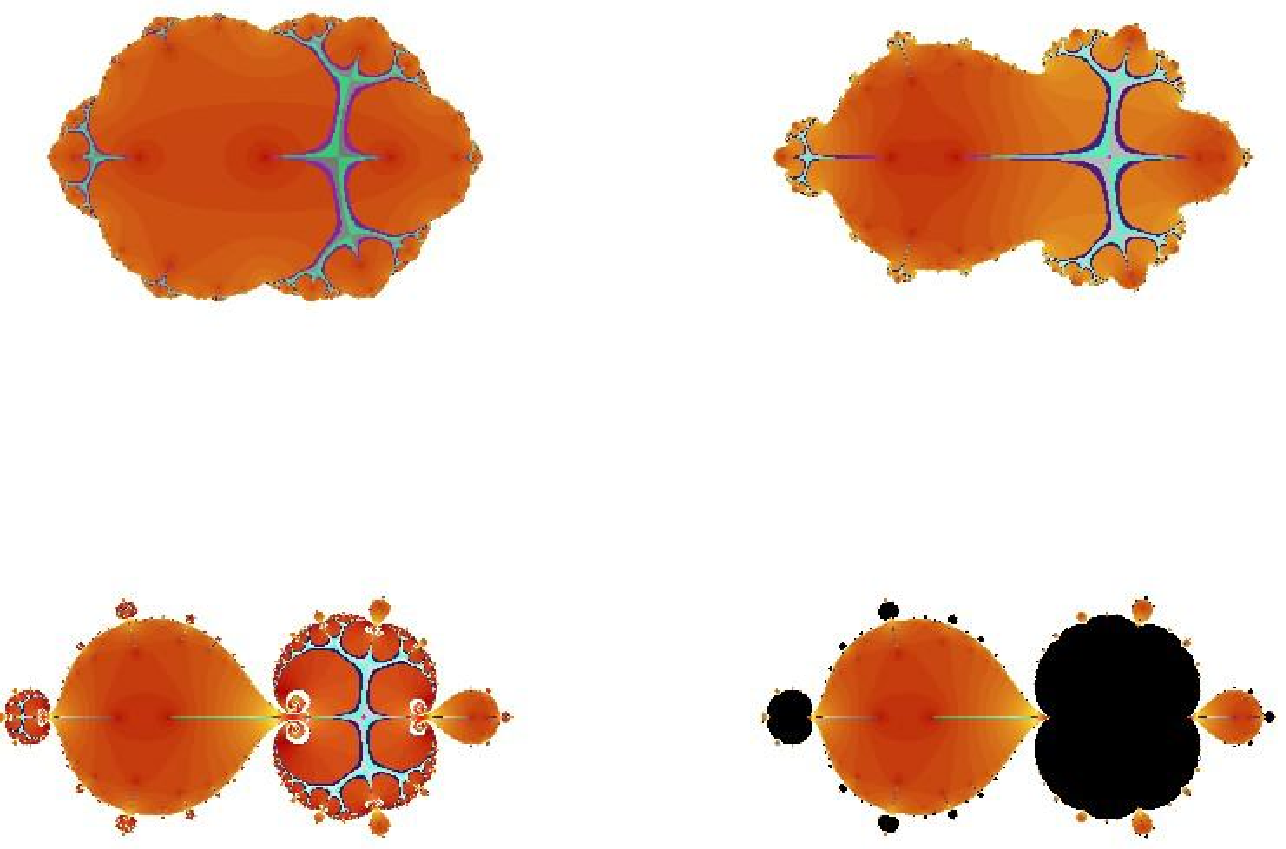}{4.0}{Evolution of spinning.  Note the apparent
collision of the pair of repelling fixed points, the creation of
the bottleneck between the spirals, and the resulting parabolic
implosion.}  
\gap

\subsection{Main results}

We will consider only the case of spinning critical points in attracting
basins.  That is, we consider spinning a single puncture
$\Psi(c)$ around a curve $\gamma \subset S^\sharp$, where $S$ is a complex
torus
$T$ with punctures.  We do not require that $c$ is in the {\em
immediate}  basin $B$ of the attractor (in this case, $c$ would not be
visible).  Recall that the torus
$T$ has a canonical homology class
$\alpha$ represented by a counterclockwise oriented simple closed curve
surrounding the attractor which is round in the linearizing coordinates.
\gap

\noindent{\bf Standing assumptions.}  In each of the results below, we
assume that a rational map $f$ is given, and the triple $(S,c,\gamma)$
for which we spin is chosen as follows:
\be
\item[{\bf A1.}] $a$ is an attracting {\em fixed} point, $B$  is its
immediate basin, $S=\hat{B}/f$ is the quotient surface,
$\Psi: \hat{B} \to S$ is  the projection, and $c$ is a
critical point whose grand orbit passes through $B$;

\item[{\bf A2.}] $\gamma \subset S^\sharp \equiv S\union\{\Psi(c)\}$ is a
simple closed curve containing $c$ such that $[\gamma] \cdot [\alpha] =
+1$, where
$[\alpha], [\gamma] \in H_1(T,\Z)$ are the corresponding classes and "
$\cdot$ " is the signed homological intersection number.

\item[{\bf A3.}]  the grand orbit of $c$ does not contain
other critical points, and there exists a critical point $b \in B$ which is
distinct from $c$ and which is visible with respect to $\gamma$.  

\eb 

The sign conventions on the intersection number mean in particular that if
we lift
$\alpha, \gamma$ under the projection $p:
\C - \{0\} \to T$ to oriented curves $\widetilde{\alpha},
\widetilde{\gamma}$ in the linearizing coordinate plane for
$a$, then $\widetilde{\alpha}$ winds counterclockwise once about the
origin and $\widetilde{\gamma}$ is an infinite ray pointing away from the
origin, invariant by $z\mapsto f'(a)z$.
\gap

\noindent{\bf Remarks:}  {\bf A1.} keeps our exposition
free of additional notation.  {\bf A2.} implies that the spun
critical points are pushed away, rather than toward, attractors.  {\bf
A3.} keeps the discussion generic and avoids the need for separate 
consideration of a plethora of combinatorially distinguished special
cases.  
\gap

\noindent{\bf Genericity.  Critical orbit relations.}   A rational map
(respectively, a polynomial) is {\em critically generic} if every
critical point (respectively, every finite critical point) is simple.   A
rational map (polynomial) has {\em no critical orbit relations} if  the
grand orbits of any two distinct (finite) critical points are disjoint,
and the forward orbit of every (finite) critical point is infinite.
\gap

\begin{thm}[Lands or diverges]
\label{thm:lands_or_diverges} Suppose $f$ is a critically generic
hyperbolic rational map with no critical orbit relations, and $(S,c,\gamma)$ satisfies the
standing assumptions.

Then the spinning ray $\sigma_\gamma$ either has a unique limit  in
$\Rat_d/\Aut(\IP^1)$, or else converges to infinity.  Furthermore, the
limit depends only on $f$ and the homotopy class of $\gamma$ in
$\pi_1(S^\sharp, c)$.
\end{thm}
(By converging to infinity, we mean that given any compact
subset $K$ of $\Rat_d/\Aut(\IP^1)$, there is an $R=R(K)$ such that
$[f_t]=\sigma_\gamma(t)  \not\in K$ whenever $t > R$.)

\begin{thm}[Lands for polynomials]
\label{thm:lands_for_polynomials} Let $f$ be a hyperbolic polynomial with
connected Julia set and no critical orbit relations, but possibly having
multiple critical point.  Then the spinning ray $\sigma_\gamma$ has a
unique limit in $\Poly^\times _d /\Aut(\C)$ which depends only on $f$ and the
homotopy class of $\gamma$ in $\pi_1(S^\sharp, c)$.  
\end{thm}

We note that Cui, by very different methods, has announced sufficient
criteria for existence of spinning limits in which {\em multiple}
critical points are spun \cite{cui:limits}.  

The following theorem states what dynamical features are inherited by
limits of spinning.   

\begin{thm}[Spinning limit inherits a large part of the dynamics]
\label{thm:closure} \mbox{}\\Suppose $F$ is an arbitrary rational map and
$(S,\gamma,c,a,B)$ satisfies the standing assumptions.
Let $F_{t_n}=H_{t_n} \circ F \circ H_{t_n}^{-1}
\in \sigma(t_n)$ be rational maps produced by spinning $c$ around
$\gamma$, where
$\{t_n\}_{n=0}^\infty$ is any sequence of real numbers.
Suppose
$F_{t_n} \to G \in \Rat_d$.
Then:
\be
\item Let $\WWW_0$ denote the union of those Fatou components which do
not iterate to $B$ or to a Siegel disk.
 Then there is a holomorphic embedding $\JJJ: \WWW_0 \to \IP^1$ such that
$\JJJ \circ F = G \circ \JJJ$. 

\item If $a_{t_n} =H_{t_n}(a)$, then the $a_{t_n}$ are attracting fixed points
of constant multiplier, and after possibly passing to a subsequence,
$a_{t_n} \to a_\infty$, an
attracting fixed point of
$G$ of the same multiplier.

\item Let $B_{t_n}$ be the immediate basin of $a_{t_n}$ under
$F_{t_n}$, and let $B_\infty$ be the immediate basin of
$a_\infty$ under $G$.

 There exists an open subset $\UUU$ contained in the grand orbit
of the basin $B$ and holomorphic embeddings $J_{t_n}: \UUU \to \IP^1$
such that: 
(i) $F(\UUU) \subset \UUU$; 
(ii) $\UUU$ contains the attractor $a$ and all critical points 
converging to $a$ except those in the grand orbit of $c$;  
(iii) $J_{t_n} \circ F = F_{t_n} \circ J_{t_n}$; 
(iv)  after passing to a subsequence, the embeddings
$J_{t_n}$ converge uniformly on compact subsets to  an embedding
$J: \UUU \to \IP^1$ satisfying $J \circ F = G \circ J$.

\eb
\end{thm}

The orientation conventions imply that the spun critical point moves away
from the attractor.  The following results describes the possibilities
for the new dynamics arising in the limit.  

\begin{thm}[Spinning limit possibilities for $\mathbf{c}$]
\label{thm:coarse}
Let $F$ be an arbitrary rational map and $(c,\gamma)$ satisfy the
standing assumptions, where $c$ is visible after $r$ steps.  Let
$F_{t_n}=H_{t_n} \circ F \circ H_{t_n}^{-1} \in \sigma(t_n)$ be
rational maps produced by spinning $c$ around $\gamma$, with $t_n \to
+\infty$ as $n \to \infty$.   Suppose $F_{t_n} \to G \in \Rat_d$
and $H_{t_n}(c)\to c_\infty$.  

Then either $G^r(c_\infty)$ lies in a fixed parabolic basin of
multiplier 1, or $G^r(c_\infty)$ is a repelling or indifferent fixed
point.   
\end{thm}

Under further assumptions, we can make more precise the connection
between visibility and the limiting dynamics.   Combined with Theorem
\ref{thm:closure}, the previous theorem shows:

\setcounter{cor}{4}
\begin{cor}[Geometrically finite limits]\label{finite}
Under the hypothesis of Theorem \ref{thm:coarse},
assume in addition that $F$ is hyperbolic.  Then 
\be
\item if $r=0$, i.e. $c$ is visible, then $c_\infty$ lies in the
immediate basin of attraction of a parabolic fixed point with
multiplier $1$; 
\item if $r>0$, i.e. $c$ is visible after $r\geq 1$ steps, then
$G^r(c_\infty)$ is a repelling fixed point of $G$.   
\eb
In particular,  the limit $G$ is
geometrically finite, possessing a single critical point $c_\infty$ which
does not converge to an attractor.
\end{cor}

The arguments used to prove the above results do not identify how the
parabolic fixed point in case 1 is created.  With more work, we have the
following.  

First, some notation.  Let $A_c$ be the central subannulus of $A$ on
which the spinning map
$h$ is holomorphic (see Figure \ref{fig:torusplot.eps}). Denote the
boundary components of $A_c$ by $\hat{\delta}^\pm$. Their lifts to the dynamical space
have a unique component $\delta^\pm$ with $a$ as one endpoint. More precisely
 $\delta^+$ joins $a$ to repelling or parabolic fixed point $u^+$;
similarly $\delta^-$ joins $a$ to a point $u^-$ (see Lemma
\ref{lemma:endpoints}). The points $u^+$ and
$u^-$ may or may not coincide.  In Figure \ref{fig:vis.eps}, the points
$u^{\pm}$ are the points in the Julia set directly above and below,
respectively, the attractor $a$.  As we perform spinning along $\gamma$,
we obtain $F_t=H_t\circ F\circ H_t$.  Denote by $\delta^\pm_t, \; a_t,
\; \ldots $, etc. the images $H_t(\delta^\pm), \; H_t(a), \; \ldots $
etc. .

\setcounter{thm}{5}
\begin{thm}[How parabolic is created]
\label{thm:fine_generic} Assume that $F$ is hyperbolic,
$(S,c,\gamma)$ satisfies the standing assumptions,  and $F_{t_n}\to G$.
Assume, by taking subsequences if necessary, 
$u^{\pm}_{t_n} \to u^{\pm}_\infty$ and $c_{t_n}\to c_\infty$.

Assume that $c$ is visible. Then
\be
\item The map $G$ has a unique parabolic  $\Omega$. It is fixed, contains $c_\infty$ and
intersects the orbits of no other critical points. 
\item The parabolic point $v$ of $\Omega$ lies on the boundary of $ B_\infty$.
\item $u^+_{\infty} = u^{-}_\infty = v$ and $u^+ \neq u^-$.
\item the multipliers 
$\lambda^\pm_{t_n}=F'_{t_n}(u^{\pm}_{t_n})$ satisfy 
\[ m < \Real\left(\frac{1}{1-\lambda^+_{t_n}} +
\frac{1}{1-\lambda^-_{t_n}}\right) < 1\]
for some real number $m<1$, and therefore $\lambda_{t_n}\to 1$
tangentially as
$n\to
\infty$.  
\eb

  If $c$ is visible after $r\geq 1$ steps then $u^+=u^-$,
$u^\pm_\infty=G^r(c_\infty)$ and $u^+_\infty\in \partial B_\infty$.
\end{thm}

The proof relies on a soft but subtle analysis of certain geometric
limits, developed in \S \ref{app:geomlims}.  The proof also shows

\begin{thm}[A case of divergence]\label{thm:diverge} Assume that $f$ is a
hyperbolic rational map (not necessarily critically generic),
$(S,c,\gamma)$ satisfies the standing assumptions,
 and $c$ is visible.
If $u^+=u^-$, then $\sigma_\gamma$ converges to infinity.
\end{thm}

\subsection{Outline of paper}  

The proof of Theorem \ref{thm:lands_or_diverges} comprises the following 
steps, which include proving Thms.
\ref{thm:closure}, \ref{thm:coarse}. Theorem
\ref{thm:lands_for_polynomials} is proved in essentially the same way.  

Step 1.  The spinning path $\sigma$
lifts to a continuous path $\sigma^\times : [0, \infty) \to X$, where $X$ is a Riemann
surface lying in a suitable space of maps $f^\times$ with normalized
and marked critical points.   We give two constructions of $X$, one in \S
\ref{secn:slice} using holomorphic motions, and a second in \S
\ref{secn:teich} using properties of puncture-forgetting maps between
Teichm\"uller spaces.  

Step 2. $\sigma(t_n)\to [g]$ in $\Rat_d/\Aut(\IP^1)$ if and only if 
$\sigma^\times (t_n)\to g^\times$ in $X$. This follows from Thm. 
\ref{thm:closure}; see \S \ref{secn:closure}.  

Step 3. The set of accumulation points of $\sigma^\times(\R)$ (with
respect to the topology of $X$)
 is discrete in $X$.  This follows immediately from Thm.
\ref{thm:coarse}, proved in \S
\ref{secn:coarse}, since multipliers vary holomorphically.

Step 4.  We assemble the results in \S
\ref{secn:main} to prove that the set of accumulation
points of $\sigma^\times([0, \infty)$  (with
respect to the topology of $X$) is either empty or one point, using in an
essential fashion the connectedness of the image of the spinning ray. In
the latter case the spinning ray converges in $X$.  Hence by Step 2, the
spinning ray  converges in
$\Rat_d/\Aut(\IP^1)$ .  

Step 5. Independence from non-canonical choices is shown 
in \S \ref{secn:teich}, using Teichm\"uller  theory and results of Bers
and Nag.

Appendix \ref{app:geomlims} develops the theory
of geometric limits of invariant strips needed for the proof of Thms.
\ref{thm:fine_generic} and \ref{thm:diverge}. Appendix \ref{app:analytical} contains miscellaneous analytical results
used in several places.   

\subsection{Acknowledgements}

We are very grateful to Curt McMullen for many useful conversations, and
to the Universit\'e de Cergy-Pontoise for financial support. 

\section{Construction of $\mathbf{X}$}
\label{secn:slice}

The main result of this section is the construction of a certain Riemann
surface $X$ consisting of rational maps with normalized marked
critical points.  The surface $X$ is defined implicitly by fixing the
dynamical behavior of all but the spun critical point.  It will contain a
lift of the spinning path, and any limit of spinning, as we shall
later show.  

Let $S, B, a, c$, etc. be as in the setup for spinning.  Let $\phi: (B,
a) \to (\C, 0)$ be a linearizing coordinate and $p: \C-\{0\} \to T$ be
projection from this coordinate to the quotient torus.  
Since the spinning homeomorphism $\widetilde{h}^t: T \to T$ acts
trivially on the homology of the quotient torus, it lifts under $p$ to a
map $\widetilde{h}_t: (\C, 0) \to (\C, 0)$.  Let $a_t = H_t(a)$ and $B_t =
H_t(B)$ be the corresponding attractor and basin for $f_t$.  

\begin{lemma}
\label{lemma:commutes} For all $t \in \R$, $\widetilde{h}_t$ 
is the identity outside of $p^{-1}(A)$, and the map
$\phi_t:=\widetilde{h}_t \circ \phi|_{\widehat{B}}
\circ H_t^{-1}: \widehat{B}_t \to \C$ is a holomorphic linearizing map
conjugating
$f_t$ to multiplication by $\lambda$.  In particular the multiplier of
$a_t$ is again
$\lambda$. 
\end{lemma}

\xincludeps{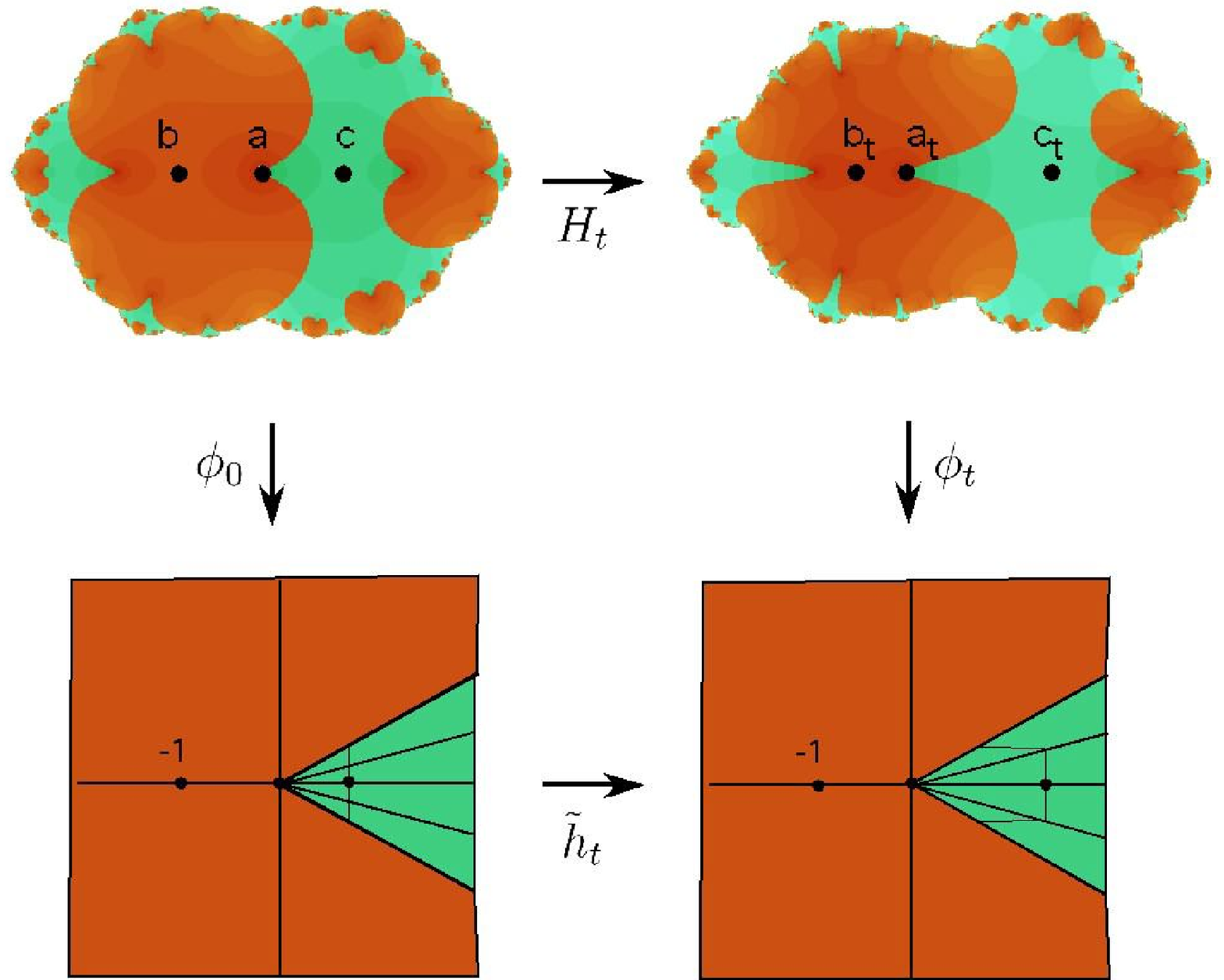}{3.5}{The qc map $\widetilde{h}_t$. The region
$p^{-1}(A)$ is the large light wedge on the right.}

\pf See Figure \ref{fig:compfig.eps}.
By construction the map $\widetilde{h}_t \circ \phi
\circ H_t^{-1}: B_t \to \C$ is well-defined, holomorphic with
respect to the standard conformal structure and  conjugates $f_t$ to
multiplication by
$\lambda$. \qed

Our strategy for creating $X$ is to work within a suitable space
$\GRatmn$  in which the critical points are marked and the maps are
normalized by e.g. conjugating so three critical points are at $0, 1,
\infty$.   Given the data defining spinning, label the critical points
of $f$ and normalize to get a map $f^\times = f^\times_0 \in \GRatmn$.
We will produce spaces 
\[ Z(f^\times_0) \supset Y(f^\times_0) \supset X(f^\times_0)\] where
$Z(f^\times_0)$ is open in $\GRatmn$ and
$X(f^\times_0), Y(f^\times _0)$ are respectively the connected components
containing $f^\times_0$ of the fibers of two holomorphic maps
$\Lambda, \Phi$ defined on $Z(f^\times_0)$ and
$Y(f^\times_0)$, respectively.
\gap

\subsection*{Marked generic rational maps.  The space
$\mbox{\bf GRat}\mathbf{^{\times,*} _d}$.}

Fix a degree
$d
\geq 3$.  Let $\GRat \subset \Rat_d$ denote the subspace of critically
generic rational maps (i.e. all critical points are simple).  
Clearly, this is open and dense.  By passing to a finite covering, we may
assume that the locations of critical points are globally defined
functions.  The Lie group $\Aut(\IP^1)$ will then act freely on this
space, and it is then simple to show that the quotient $\GRatmn$ is a
complex manifold which may be concretely realized as a subset of
$\Rat_d \times \left(\IP^1\right)^{2d-2-3}$.  For the remainder of
this section, we will work exclusively within the space
$\GRatmn$.

\gap

\subsection*{Marking attractors.  The space
$\mathbf{Z(f^\times_0)}$.}

Here, we describe a general construction
of analytic subsets of parameter spaces in which the behaviors of some, but
not all, critical points are held fixed.
\gap

\noindent{\bf Data for the definition.}   Write
\[ C \equiv \{1, 2, 3, ..., 2d-2\}\]
and choose a decomposition
\[ C =I \sqcup J \sqcup K, \;\;\;I \neq \emptyset.\]
Choose a function $\omega: J \to I$ in case $J$ is nonempty. {\em
Throughout the remainder of this section, we assume that these choices
have been given.}
\gap

\noindent{\bf Definition of $Z$.}   Let $Z=Z(I,J,K, \omega)$ denote the
following subspace of $\GRatmn$:  a normalized, marked map $f^\times$
belongs to  $Z$ if and only if
\bi
\item for $i \in I$, the critical point $c_i$ is in the {\em immediate}
basin $B_i$ of an attracting, but not superattracting, cycle
$\genby{a_i}$, $c_i$ has infinite forward orbit, and for
$i \neq i' \in I$ the basins of $\genby{a_i}$ and
$\genby{a_{i'}}$ are disjoint;
\item for $j \in J$, the critical point $c_j$ is in the basin (not
necessarily the immediate basin)  of
$\genby{a_{\omega(j)}}$, where $\omega: J \to I$ the given function; the
forward orbit of $c_j$ is infinite, for $j \neq j'$ the grand orbits
of $c_j$ and $c_{j'}$ are distinct, and the grand orbits of $c_j$ and $c_{\omega(j)}$
are disjoint;
\item for $k \in K$, there are no restrictions on the behavior of $c_k$.
\ib

Note that a given $f^\times_0 \in \GRatmn$ for which the underlying map
$f_0$ is hyperbolic and without critical orbit relations  can be regarded as an element of $Z$ for a
variety of different choices of subsets $I,J$, and that such a choice
determines the function
$\omega:J \to I$ uniquely.  Also, the space $Z$ contains many
combinatorially distinguished connected components, since we have not
specified e.g. how many iterates are needed for $c_j$ to map into the
basin of $a_{\omega(j)}$.

We let $Z(f^\times_0)$ denote the connected component of $Z$
containing $f^\times _0$.  The space $Z$ is open in
$\GRatmn$, and is therefore a complex manifold of dimension $2d-2$.  In
particular, $Z(f^\times_0)$ is a complex manifold of dimension $2d-2$.

\subsection*{Fixing the multipliers.  The space $Y(f^\times_0)$.}   For
$f^\times \in Z$ and $i
\in I$, let
$a_i(f^\times) \in \IP^1$ denote the location of the attracting periodic
point whose immediate basin contains
$c_i$, and let
$p_i(f^\times)$ be the period of this attractor. Clearly, these are
functions of $f^\times$. Hence the locations of each point in the
attractors  $\genby{a_i}$ , and the multiplier
$\lambda_i$ of this attractor, are in fact a function of
$f^\times$.

We denote by $\Lambda$ the multiplier map:
\[ \Lambda: Z \to (\Delta^*)^{I}\ ,\quad f^\times \mapsto \left(\lambda_i(f^\times)\right)_{i \in I}.\]

\begin{prop}
\label{prop:section_for_lambda} The map $\Lambda$ is holomorphic and
admits local holomorphic sections.  That is, given any point $\vec{\lambda}_0=
(\lambda_i)_{i
\in I} \in (\Delta^*)^I$, and any map $f^\times$ with
$\Lambda (f^\times) =
\vec{\lambda}_0$, there is a neighborhood $U$ of $\vec{\lambda}_0$ and a
holomorphic map
\[ \Sigma: U \to Z \] such that
\[ \Lambda \circ \Sigma = \id_U.\]
\end{prop}

 From the theory of holomorphic functions of several complex
variables, we immediately obtain (see Corollary C.10, p. 23 of
\cite{gunning:scv:I}).

\begin{cor}
\label{cor:lambda_fiber} Given any $f^\times _0 \in Z$, the fiber
\[ \Lambda ^{-1}(\Lambda(f^\times _0))\] is a complex manifold of dimension
$2d-2-|I|$ which is closed subset of $Z$, with respect to the induced topology to $Z$.
\end{cor}

Note that the fibers $\Lambda^{-1}(\Lambda(f^\times _0))$ need not be
closed in $\GRatmn$.

Given $f^\times _0 \in Z$, we let
$Y(f^\times_0)$ denote the connected component of the fiber
$\Lambda^{-1}(\Lambda(f^\times _0))$ which contains
$f^\times _0$.  This is again closed in $Z$.
\gap

\pf {\bf (of Prop. \ref{prop:section_for_lambda}):}

Denote by $\IH=\{x+iy,x>0\}$ the right half plane. For any $s\in \IH$, 
define  a homeomorphism $l_s :  \C\to \C$ by
$$l_s(z)=l_s(re^{2\pi i \theta})=r^se^{2\pi i \theta}=z\cdot r^{s-1}=z\cdot e^{(s-1)\log r}\ .
$$
It is a quasi-conformal map on $z$, and it depends holomorphically on $s$.
An easy calculation shows that $l_s(\lambda z)=\lambda|\lambda |^{s-1}l_s(z)$, for any $\lambda\ne 0$.

Suppose $f^\times$ is now given. For $s=(s_i)\in \IH^I$, set
$\vec{\lambda}(s)= (\lambda_i|\lambda_i|^{s_i-1})_{i\in I}$.
Fix
temporarily $i \in I$. Let $a_i \in \IP^1$ be
the attractor whose immediate basin $B_i$ contains the critical point
$c_i$.  Let $\widetilde{B}_i$ be the entire basin of $a_i$. Choose a holomorphic map
\[ \psi_i: (B_i, a_i) \to (\C,0)\] satisfying
\[ \psi_i \circ f^{\circ p_i}(z) = \lambda_i \psi_i(z)\ .\]
Entend it then to  $\widetilde{B}_i$ by the following recipe: 
\[ \psi(z) =
\lambda^{-[\frac{n}p]}\psi(f^{\circ n}(z)) \] where $n$ is any
nonnegative integer for which $f^{\circ n}(z)
\in B'$ and $[n/p]$ is the greatest integer less than or equal to
$n/p$. The extended $\psi$ satisfies the same functional equation and maps
a grand orbit of $f$ onto a grand orbit of $\lambda z$. 

For any $s=(s_i)\in \IH^I$, we define a new complex structure $\sigma(s)$ as follows: for each
 $i\in I$, $\sigma(s)|_{\widetilde{B}_i}$ is the pull-back of the standard complex
structure under $l_{s_i}
\circ \psi_i$ , and  $\sigma(s)|_{\cbar-\bigcup_i\widetilde{B}_i}$ is the
standard complex structure. This complex structure is $f$-invariant by
construction. Denote by $h_s$ the integrating map fixing $0,1,\infty$, 
by $g_s$ the rational map $h_sfh_s^{-1}$, and by $g_s^\times$
the marking $(g_s,h_s(\vec{c}(f^\times)))$. See the following diagram:
$$\begin{array}{cclcl}
f,\ \cbar\supset \bigcup_i \widetilde{B}_i \ni a_i &  \vspace{0.3cm}  \stackrel{h_s}{\longrightarrow} & 
 \cbar\supset \bigcup_i \widetilde{B}_i(s)\ni a_i(s),\ g_s \\  
\downarrow_{\prod_i \psi_i}  && \downarrow_{\prod_i(l_{s_i}\psi_i) h_s^{-1}} \vspace{0.3cm}   \\
\C^I\vspace{0.3cm},\  \prod_i(\lambda_i z) & \stackrel{\prod_i l_{s_i}}{\longrightarrow} & 
 \C^I,\ \prod_i(\lambda_{s_i}z)  \end{array}\ .$$

Note that for $s=(1,\cdots,1)$, $l_{s_i}=id$, $h_s=id$ and $g_s=f$.
Otherwise $s\mapsto g_s^\times$ is holomorphic, with $\Lambda(g^\times_s)=\vec{\lambda}(s)$.
Note  that the map 
$w: s\mapsto \vec{\lambda}(s)$ is locally bi-holomorphic mapping a small neighborhood
of $(1,\cdots,1)$ onto a neighborhood $U$ of $\vec{\lambda}_0$. Moreover the maps $g_s^\times$ are in
$Z$. As a consequence, $\Sigma: U\to Z$, $\Sigma(\vec{\lambda})=g^\times_{\omega^{-1}(\vec{\lambda})}$
 is a holomorphic section of $\Lambda$.
\qed

\gap
\subsection*{Fixing critical points in linearized coordinates.\\  The
space $\mathbf{X(f^\times_0)}$.}  

On the space $Z$, there is another map defined as follows.  Fix
temporarily $i \in I$ and $f^\times \in Z$.  Let $a_i \in \IP^1$ be
the point in the attracting cycle $\genby{a_i}$  whose immediate basin $B_i$ contains the critical point
$c_i$.   There is a unique normalized  holomorphic map
\[ \psi_i: (B_i, a_i) \to (\C,0)\] satisfying
\[ \psi_i \circ f^{\circ p_i}(z) = \lambda_i \psi_i(z),\quad \mbox{and } \psi_i'(a_i)=1\ .\]
   By hypothesis, $c_i$ has infinite forward orbit,
so $\psi_i(c_i) \neq 0$.  Rescale to set 
\[ \phi_i(z)=-\frac{1}{\psi_i(c_i)}\psi_i(z).\]
 Then
\[ \phi_i: (B_i, a_i, c_i) \to (\C, 0, -1).\]
It is a linearizing map. We then extend it to $\widetilde{B}_i$ as in the proof of 
Proposition \ref{prop:section_for_lambda}.
Recall that $\phi_i$ maps grand orbits of $f$ onto 
$\lambda_i$-orbits.

From the definition of the subspace $Z$, recall that 
\bi
\item if $\omega(j)=i$, then the critical point $c_j$ lies
in $\widetilde{B}_i$.  Hence the value $\phi_i(c_j)$ makes
sense.  
\item the critical points $c_j$ have infinite forward orbits.  Hence if
$\omega(j)=i$, then  $\phi_i(c_j) \neq 0$.
\item for $j \neq j'$, the critical points $c_j$ and $c_{j'}$ have
distinct grand orbits.  Suppose in addition that $j \neq j'$ and
$\omega(j)=\omega(j')=i$.  Then $\phi_i(c_j), \phi_i(c_{j'})$ 
have distinct $\lambda_i$-orbits.  In particular, $\phi_i(c_j) \neq
\phi_i(c_{j'})$.
\item the critical points $c_i$ and $c_j$ have distinct grand orbits.
Hence if $\omega(j)=i$ then $\phi_i(c_j)\neq -1$, and $\phi_i(c_j)$ and
$-1$ have distinct 
$\lambda_i$-orbits.  
\ib
\gap

\noindent{\bf Notation.}  We set \[ \mathbf{C}^i = (\C^*-\{\lambda_i^n(-1),n\in \Z\})^{|\omega^{-1}(i)|} 
- \mbox{big diagonal}\]
where the big diagonal is the locus where two or more coordinates have the same $\lambda_i$-orbit.  
\gap

Thus for each $i \in I$ we have a function 
\[ \Phi_i: Z \to \mathbf{C}^i\quad
\mbox{given by}
\ \Phi_i(f^\times) = \left(\phi_i(c_j(f^\times))\right)_{j \in \omega^{-1}(i)}.\]
Putting these together, we have a function 
\[ \Phi : Z \to \prod_i \mathbf{C}^i\subset \C^J\quad
\mbox{given by}
\ \Phi(f^\times) = (\phi_i(c_j(f^\times)))_{j \in J,\,i=\omega(j)}\]
which records the locations of the critical points $c_j$ in the
linearizing coordinates.

\begin{prop}
\label{prop:section_for_phi} The map $\Phi$ is holomorphic.  The
restriction of
$\Phi$ to any fiber of $\Lambda$ admits local holomorphic sections.  That
is, given any point $\ww_0 \in \prod_i \mathbf{C}^i$, and any map
$f^\times $ with
$\Phi(f^\times) = \ww_0$, there is a neighborhood
$U$ of $\ww_0$ and a holomorhic map
$ \Sigma:U \to Y(f^\times)$ such that
$ \Phi\circ \Sigma = \id_U.$
\end{prop}

\begin{cor}
\label{cor:phi_fiber} Given $f^\times _0 \in  Z$, the fiber of the
restriction
\[
\left(\Phi\left|_{Y(f^\times_0)}\right.\right)^{-1}(\Phi(f^\times _0))\]
is a complex manifold of dimension $2d-2-|I|-|J|$ which is closed as a
subset of $Y(f^\times_0)$, therefore closed in $Z$.
\end{cor}

\pf {\bf (of Prop. \ref{prop:section_for_phi}):} Let $\ww_0=(w_j)_{j\in J}$.
For each  $i\in I$ and 
 each $j\in \omega^{-1}(i)$, choose a small round disk $D_j$ centered at $w_j$
 so that $\lambda_i^nD_j$ are disjoint
for distinct $n\in \Z$, and the
$\lambda_i$-orbits of $D_j$, of $D_{j'}$ and of $-1$ are mutually disjoint, for any $j'\ne j$ and 
$\omega(j')=i$.
Set $\Delta_i=\prod_{j,\omega(j)=i}\frac12 D_j$, and $U=\prod_i\Delta_i=\prod_{j\in J}\frac12 D_j$.

Fix $i\in I$. For each $j\in \omega^{-1}(i)$, define a holomorphic motion
 $M_{ij}:\frac12 D_j\times \overline{D}_j\to \overline{D}_j$ as follows:

$M_{ij}(u,w_j)=u$, $M_{ij}(u,z)=z$ for any $(u,z)\in \frac12 D_j\times \partial D_j$,
$M_{ij}$ is holomorphic on $u$ and injective on $z$, and $M_{ij}(w_j,\cdot)$ is the identity.

Let $M_i: \Delta_i\times \C\to \C$ be the following  holomorphic motion: 

\noindent$\bullet$ For each
 $j\in \omega^{-1}(i)$ and any pair $(\uu,z)\in \Delta_i\times \overline{D}_j$ with $\uu=(u_j)$,
 $M_i(\uu,z)=M_{ij}(u_j,z)$.
In particular $M_i(\uu,w_j)=u_j$.

\noindent$\bullet$  $M_i(\uu,z)=\lambda_i^n M_i(\uu,z/\lambda_i^n)$  for any $n\in \Z$ and any pair 
$(\uu,z)\in \Delta_i\times \lambda_i^n\overline{D}_j$.

\noindent$\bullet$ $M_i(\uu,z)=z$, for $ z\in \bigcup_j\partial D_j$ and for $z$  outside of the
 $\lambda_i$-orbits of $\bigcup_jD_j$.

Note that by construction $M_i(\uu,\cdot)$ commutes with the multiplication by $\lambda_i$.

Do this for every $i\in I$.

Fix now $\uu\in U$. We define a new complex structure $\sigma(\uu)$ as follows: for each $i\in I$,
$\sigma(\uu)|_{\widetilde{B}_i}$ is the pull-back by $M_i(\uu|_{\Delta_i},\cdot)$ and then by
$\phi_i$ of the standard structure, $\sigma(\uu)|_{\cbar-\bigcup_i\widetilde{B}_i}$ is the
standard structure. Such structure is $f$-invariant by construction, and is holomorphic on $\uu$.
Let $h_\uu$ be the unique integrating map fixing $0,1,\infty$. Let $g_\uu=h_\uu fh_\uu^{-1}$
and $g_\uu^\times =(g_\uu,h_\uu(\vec{c}(f^\times)))$.
Then $\Lambda(g_\uu^\times )\equiv \Lambda(f^\times)$, and $\Phi(g_\uu^\times)=\uu$. Therefore 
$\Sigma(\uu)=g_\uu^\times$ is a holomorphic section of $\Phi|_{Y(f^\times)}$.
\qed

Note that the same proof can be adapted to get similar results in a parabolic basin or a rotation domain.

Given $f^\times_0$, we let $X(f^\times_0)$ denote the connected
component of the fiber in the above corollary containing $f^\times_0$, which is again closed in $Z$.

\begin{cor}[One-dimensional]
\label{cor:one_dimensional} If $|K|=1$, i.e. if there is a single free
critical point, then
$X(f^\times_0)$ is one-dimensional.
\end{cor}

\subsubsection*{Polynomial case}
We briefly sketch the construction of analogous spaces for polynomials, having 
possibly multiple critical points.

A polynomial of degree $d$ has $d-1$ critical points, counted with
multiplicity.  Given a partition $\DDD$ of $d-1$:
\[ d-1=d_1 + d_2 + ... + d_M \]
define
\[ \Poly^\times _d (\DDD) = \left\{ f^\times(z) \equiv d\int_0 ^z
\prod_{m=1}^M (\zeta-c_m)^{d_m}d\zeta\; \right\} \leftrightarrow \{(c_m)\in\C^M\}.\]
Note that elements of $\Poly^\times_d$ are polynomials which
fix the origin, are monic, and whose critical points are labelled.  
The projection $\Poly^\times _d (\DDD) \to \Poly^\times _d /\Aut(\C)$
is thus finite-to-one.  We then use this space in place of $\GRatmn$,
and proceed to define $Z, Y, X$ as before.

The {\em connectedness locus} is the subspace of $\Poly^\times _d (\DDD)$
consisting of maps whose Julia set is connected.  Equivalently, the
orbit of every critical point $c_m$ is bounded.  Later, we will need the
following result:

\begin{lemma}
\label{lemma:conn_locus_bounded}
The connectedness locus is a bounded subset of $\Poly^\times _d (\DDD)$.
\end{lemma}

\pf
Suppose $f^\times \in \Poly^\times _d (\DDD)$ has connected Julia set.  Let
$K_f$ be the filled-in Julia set of $f$.  By a theorem of B\"ottcher, there is a
unique Riemann map
\[ (\IP^1-\cl{\Delta}, \infty) \to (\IP^1-K_f, \infty)\]
which is tangent to the identity at infinity and which conjugates $w
\mapsto w^d$ to $f$. Note that $0\in K_f$.  By the Koebe $\frac14$-theorem
(applied in the $1/z$ coordinates), the image of $\Sigma$
contains a spherical disk centered at infinity whose radius is
independent of $f$ ($\{|z|>4\}$).  Thus $K_f$ is contained in 
$\{|z|\le 4\}$, independent of $f$.  Since the critical points $c_m$ of
$f$ are contained in $K_f$, the lemma follows.
\qed

\subsection*{Application to spinning}
Let $f_0\in \GRat$, $\gamma$ be as in the setup for spinning, and
$\sigma$ be the
 corresponding spinning path.
 Assume that  $f^\times_0 \in \GRatmn$ is a representative of $f_0$. In other words,
we choose a labelling $c_1,\cdots,c_{2d-2}$ of the critical points of $f_0$, and we normalize
$f_0$ so that $c_1=0$, $c_2=1$ and $c_3=\infty$.

\begin{prop}[Spinning path lifts to $\GRatmn$]
\label{prop:lifts_to_G} The spinning path
\[ \sigma: (\R,0) \to (\Rat_d/\Aut(\IP^1), [f_0]) \] is continuous and is the projection
of a continuous path 
\[ \sigma^\times : (\R, 0) \to (\GRatmn, f^\times_0).\]
\end{prop}

\pf Choose the quasiconformal conjugacies $H_t$
to fix $0,1,\infty$. Then $H_t$, as well as  $H_t f_0 H_t ^{-1}$, depend continuously
(actually, real analytically) on $t$.  Set
\[ \sigma^\times (t) =f^\times_t \equiv  (H_t f_0 H_t ^{-1},
H_t(\vec{c})).\] By construction, the image of
$\sigma^\times $ lies in $\GRatmn$ and projects to $\sigma$. So $\sigma$ 
is itself continuous.
\qed

Assume further that we have written the set of indices of critical points
of $f^\times_0$ as $I \sqcup J \sqcup K$ as above, such that $f^\times_0 \in Z$
and that the spun critical point $c$ has its index in $K$
(this places restrictions on the behavior of the critical points).
This then determines
the function $\omega: J \to I$.

\begin{prop}[Spinning path lifts to $\mathbf{X(f^\times_0)}$]
\label{prop:lifts_to_X}
Suppose
\[ f^\times_0 \in Z=Z(I,J,K,\omega).\]  Then the lift
$\sigma^\times $ of the spinning path lies in $X(f^\times_0)$.
\end{prop}

\pf 
Let $i\in I$ and $c_i$ denote the $i$th critical point of $f_0^\times$.  Then
$c_i$ is in the immediate basin $B_i$ of the attractor $a_i$.

By conjugacy $H_t(c_i)$ is in the immediate basin of $H_t(a_i)$ and any $i\in I$, and $H_t(c_j)$ is in the
basin of $H_t(a_{\omega(j)})$ for any $j\in J$. So $\sigma^\times \subset Z(f^\times_0)$.

By Lemma \ref{lemma:commutes}, the multiplier of $H_t(a_i)$ is independent
of
$t$. Thus the multiplier map $\Lambda$ is constant on the spinning path
$\sigma^\times $, and so $\sigma^\times  \subset Y(f^\times_0)$.  

Recall that
\[ \phi_i: (B_i, a_i, c_i) \to (\C,0, -1)\]
is the normalized linearizing map on $B_i$.

Assume that $i\in I$ and $a_i$ does not attract the spun critical point $c$. Then for 
$ \phi_{i,t}=\phi_i\circ H_t^{-1}$, we have $\phi_{i,t}(c_{j,t}) \equiv \phi_i(c_j)$ 
,
in particular 
 $\phi_{i,t}(c_{i,t})\equiv -1$. So $\phi_{i,t}$ coincides with the normalized lineariser in the definition
of $\Phi$.

Assume now $i\in I$ such that $a=a_i$ does attract the spun critical point $c$.
Then for
$\phi_{i,t}=\widetilde{h}_t\circ \phi_i\circ H_t^{-1}$, we have $  \phi_{i,t}(c_{j,t}) 
 = \widetilde{h}_t(\phi_i(c_j)) =
\phi_i(c_j)$,
by Lemma \ref{lemma:commutes}, and the fact that $c_{j,t}$ is not in the
grand orbit of
$c$. In particular $\phi_{i,t}(c_{i,t})\equiv -1$. So, again,
 $\phi_{i,t}$ coincides with the normalized lineariser in the definition
of $\Phi$.

It follows that the function $\Phi$ is constant on $\sigma^\times $. Thus
\[ \sigma^\times  \subset X(f^\times_0)\]
and the Proposition is proved.
\qed

\noindent{\bf Remark:}
Let $f^\times_0, \gamma, \tilde{\sigma}$ be as in the example in \S
1.4 but now let $t \to -\infty$.  Using the symmetry $z \mapsto -z$ of
$f_0$ it is straightforward to verify that the resulting path is the same
as the path defined by spinning the other critical point $b$ outward
along a curve which is the image of the negative real axis under
projection from the linearizing coordinate.  In the limit, the critical
point $b$ lies in a parabolic basin (by Theorem \ref{thm:coarse}) and so
the locus $X(f^\times_0)$ of points for which $b$ converges to the origin
is indeed not closed.

\section{Limits of spinning, I} 
\label{secn:closure} 

In this section, we prove Theorem \ref{thm:closure}, which explains
what dynamical features are preserved when passing to a limit of
spinning.  
\gap

{\noindent\bf Proof}.
Conclusion (1) follows by Lemma \ref{lemma:persistence} below.
Conclusion (2) follows by Lemma \ref{lemma:commutes}.  We now prove (3).

Let $\phi:(B,a,b) \to (\C, 0, -1)$ be the normalized linearizing map for the
attractor $a$.  Extend $\phi$ to the grand orbit of $B$.  Let $p: \C^*
\to T$ be projection onto the quotient torus (i.e. identifying $z$ to $\lambda z$).  In the linearizing
coordinate plane $\C$, the set $p^{-1}(A)$ is a single, thickened
logarithmic spiral emanating away from the origin (due to the 
standing assumption A2 on $\gamma$).  The map 
\[ p^{-1} \circ p_A: (\{x+iy: |y| \leq 2l\}, 0) \to (\C^*, \phi(c))\]
conjugates translation by $-1$ to multiplication by $\lambda=F'(a)$. 
Consider the domain in the linearizing coordinate plane given by 
\[ V_0=V = \C - \cl{p^{-1}\circ p_A(\{ x+iy : x \geq -1\})}.\]
See Figure \ref{fig:holemb.eps}.  Note that $V$ is open, contains the
origin, is forward-invariant under multiplication by $\lambda$, omits
$\phi(c)$ and
$\phi(F(c))$, and contains the images under $\phi$ of all other critical
points in the grand orbit of $B$.  Let $\UUU = \phi^{-1}(V)$.  Then
conditions (i)-(ii) of the conclusion (3) hold.  

To prove (iii), let $U_0$ be the component of $\UUU$               containing $a$.  Then
$U_0$ is forward-invariant.  Let $\phi_t, \widetilde{h}_t$ be as in the
proof of Proposition \ref{prop:lifts_to_X}, let
$V_t=\widetilde{h}_t(V_0)$, and let
$U_t=H_t(U_0)$; see Figure \ref{fig:holemb.eps}.

\xincludeps{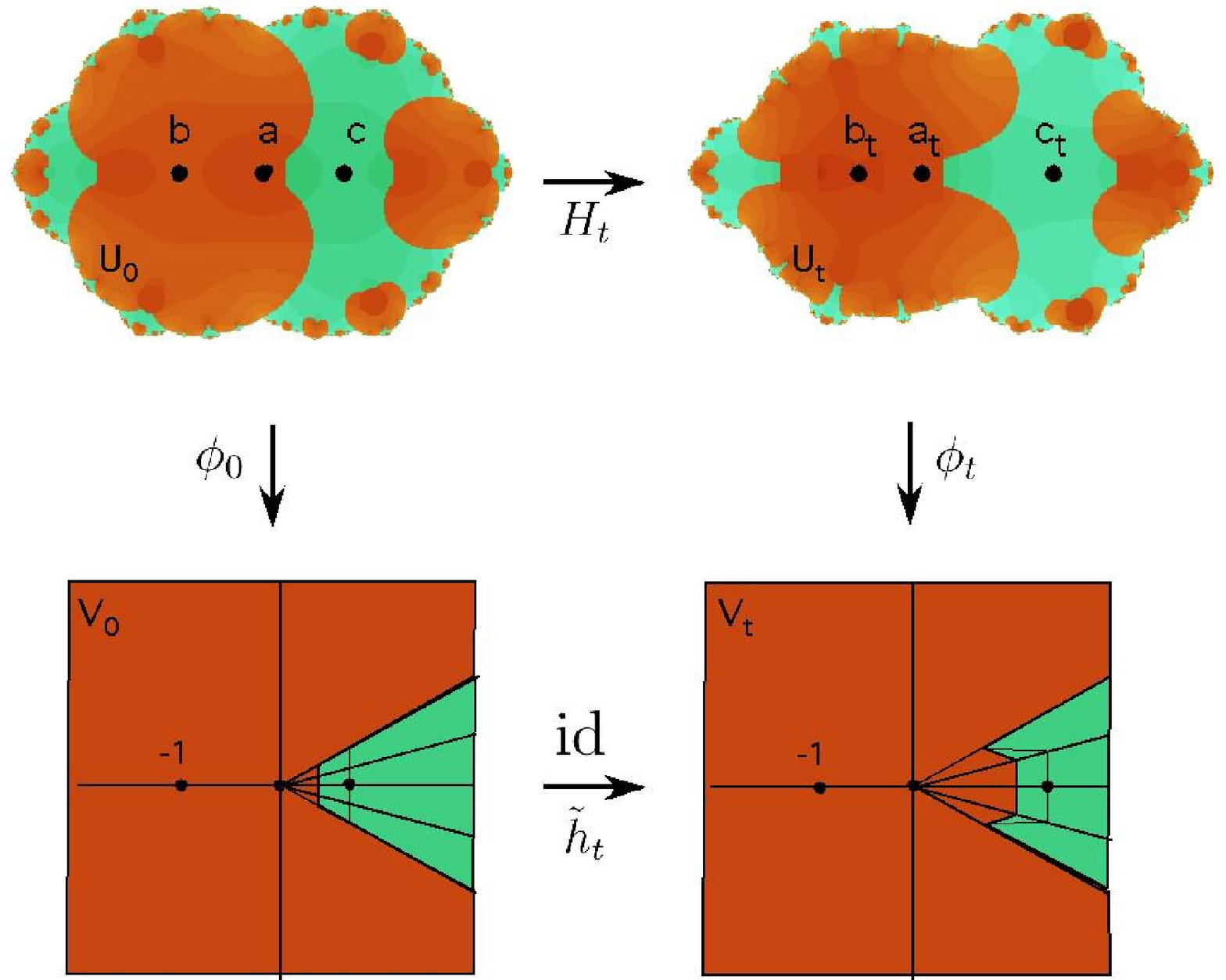}{3.5}{The identity map $\id|V_0: V_0
\hookrightarrow V_t$ is homotopic to $\widetilde{h}_t|V_0: V_0 \to V_t$
and therefore lifts.}

Fix $t \in \R$.  If $s<t$, then $V_s \subset V_t$.  Thus, the family of
maps 
\[ \widetilde{h}_s: V_0 \to V_t, \;\;\;0 \leq s \leq t\]
provides an isotopy from the restriction $\id|_{V_0}$ of the identity map
to the map $\widetilde{h}_t: V_0 \to V_t$.  
Let $B'$ be the complement in $B$ of the critical points of $F$ and
all of their backward orbits; see Lemma \ref{lemma:koenigs_cover}.  Let
$U_0' = U_0 \intersect B'$, and put $V_0' = \phi_0(U_0')$.  Then by Lemma
\ref{lemma:koenigs_cover} the restriction $\phi_0: U_0' \to V_0'$ is a
covering map.  Similarly, with the corresponding notation, $\phi_t: U_t'
\to V'_t$ is also a covering map.

For each $s$ in $0 \leq s \leq t$, the map $\widetilde{h}_s$ is the identity
on $\C-p^{-1}(A)$.  In particular, it is the identity on each puncture of
$V_0'$.   Since
$\widetilde{h}_t: V'_0 \to V'_t$ lifts under
$\phi_0, \phi_t$ to $H_t: U'_0 \to U'_t$, by lifting of isotopies we have
that $\id|V_0: V'_0 \to V'_t$ lifts to a holomorphic embedding
\[ J_t: U'_0 \to U'_t\]
such that $\phi_t \circ J_t = \id \circ \phi$ on $U$.  Since also $\phi
\circ F = \lambda \phi$ and $\phi_t \circ F_t = \lambda \phi_t$, this
implies $J_t \circ F = F_t \circ J_t$.  

Using similar reasoning, one can inductively extend $J_t$ to each
component of $\UUU$ to obtain an embedding $J_t: \UUU \to \IP^1$ such
that $J_t \circ F = F_t \circ J_t$.  Note that this implies that $J_t(b)$
is a critical point of $F_t$.  

To prove (iv), apply Lemma \ref{lemma:limits_of_models} with $W=U_0$ to conclude
that after passing to subsequences, the holomorphic maps $J_t|U_0$ converge
locally uniformly to a map $J$ which is either an embedding, or else is
a constant map with value in $\Fix(G)$.  Since $U_0$ contains both the
visible critical point $b$ and the fixed point $a$, if $J$ were constant, then $J(b)=\lim_{n \to
\infty} J_{t_n}(b)$ would be a fixed critical point of $G$, i.e. a fixed
point of multiplier zero.  On the other hand, $J(b)=J(a)=\lim_{n \to
\infty} J_{t_n}(a) = a_{t_n} = a_\infty$ is a fixed point of multiplier
$\lambda \neq 0$ by (2). This is not possible.  Hence $J$ is an embedding of $U_0$ into
$\IP^1$.  

The maps $J_{t_n}$ are actually defined on all of $\UUU$; as before we
may assume that $J_{t_n}: \UUU \to \IP^1$ converges locally uniformly to
a map $J: \UUU \to \IP^1$ satisfying $J \circ F = G \circ J$.  On each
component $W$ of $\UUU$, the map $J$ is either an embedding, or is 
constant.  Suppose $W \subset F^{-k}(U_0)$ is a component of $\UUU$.  If
$J|W$ is constant, then $G^k \circ J = J \circ F^k$ is constant on $W$. 
Since $F^k$ is open, and $F^k(W) \subset U_0$, this implies $J|U_0$ is
constant, a contradiction.  
\qed

\section{Limits of spinning, II}
\label{secn:coarse}

Here, we study what new dynamical features develop in limits of spinning.
These are summarized in Theorem \ref{thm:coarse}, which we now prove.  

\pf Theorem \ref{thm:closure} implies that after
passing to subsequences, $a_{t_n}, b_{t_n}, c_{t_n} \to a_\infty,
b_\infty, c_\infty$ where $a_\infty, b_\infty, c_\infty$ are 
distinct.  Let $M_{t_n} \in \Aut(\IP^1)$ be the unique map sending
$a_{t_n}, b_{t_n}, c_{t_n} \mapsto 0, -1, +1$.  Then $M_{t_n} \to
M_\infty$, and so
$M_{t_n} \circ F_{t_n} \circ M_{t_n}^{-1} \to M \circ G \circ
M^{-1}$.  Hence we may assume $a_{t_n}, b_{t_n}, c_{t_n} = 0, -1, +1$.

Let $\phi: (B,0,b=-1) \to (\C, 0, -1)$ be the unique normalized
linearization map for the basin $B$.  Recall that $\Psi = p \circ \phi: B
\to B/F = T$, a complex torus, and that $A$ is an annulus on $T$ with
core $\gamma$.  The statement below makes precise the idea that, under
these assumptions,
  points
which get ``spun'' move off to infinity from the point of view of the attractor
at the origin.

\begin{prop}[Spun points tend to infinity]
\label{prop:to_infty}
Let $z \in \Psi^{-1}(A)$, and suppose $z_\infty$ is any limit point of
$z_{t_n} = H_{t_n}(z)$.  Then $z_\infty$ is not in $B_\infty$.
\end{prop}

\pf By Lemma
\ref{lemma:commutes}
\[ \phi_{t_n}(z_{t_n}) = \phi_{t_n}(H_{t_n}(z)) =
\widetilde{h}_{t_n}(\phi(z)) \to \infty. \]
For example, if $\Psi(z)$ lies in the central subannulus of $A$ where
the conjugacies $h_{t_n}$
are conformal then
$\widetilde{h}_{t_n}(\phi(z))=\lambda^{-t_n}\phi(z)$ which tends to
infinity
as
$t_n
\to +\infty$.  Here, we have made nontrivial use of the hypotheses
that $\gamma$ is
oriented outward, and $t_n \to +\infty$.

Let $J: \UUU \to \IP^1$ be the partial conjugacy from $F$ to $G$ given by
Theorem \ref{thm:closure}.  We use now the assumption A3 that
$B$ contains a visible critical point $b$. Since $b \in \UUU$,
$J(b)$ makes sense.  Let $b_\infty = J(b)=\lim_{n
\to \infty}J_{t_n}(b) \in B_\infty$.  Since $b$ is visible, it has an
infinite forward orbit, hence so does $b_\infty=J(b)$.
If $\psi_{t_n}, \psi_\infty$ denote the linearizing maps normalized
to have derivative one at
the origin, then by Lemma
\ref{lemma:holo_koenigs}, $\psi_{t_n} \to
\psi_\infty$ uniformly on compact subsets of $B_\infty$.  In particular
\[ \psi_{t_n}(b_{t_n}) \to \psi_\infty(b_\infty) \neq 0\]
since $b_\infty$ has infinite forward orbit.
Hence as $n \to \infty$, the maps
\[ \phi_{t_n}(z)= -\frac{1}{\psi_{t_n}(b_{t_n})}\psi_{t_n}(z)\]
converge as well.  Thus $z_\infty \in B_\infty$ would imply that
$\phi_{t_n}(z_{t_n})$ converges
to a finite value, which by the previous paragraph is impossible.
\qed

\noindent{\bf Proof of  Theorem \ref{thm:coarse}}.  Let $A_c$ be the
central subannulus of $A$ on which the spinning map
$h$ is holomorphic (see Figure \ref{fig:torusplot.eps}).  Let $W$ be the
component  of
$\Psi^{-1}(A_c)$ whose closure contains the attractor $a$. By
the standing assumption (A2)  $W$ exists, is unique, and $F(W)\subset
W$.  (In Figure
\ref{fig:compfig.eps}, the subset $W$ is a slightly skinnier version of
the prominent light region on the right-hand side of the upper left
image, and $c$ is visible.) Suppose $c$ is visible after $r\geq 0$ steps.  Then $r$ is the
smallest nonnegative integer for which $F^r(c) \in W$.    The
conjugacies $H_{t_n}$ are holomorphic on $W$.  Lemma
\ref{lemma:limits_of_models} implies that after passing to a
subsequence, the maps $H_{t_n}|_W$ converge to a limit $H_\infty$ which is
either univalent and conjugates $F|W$ to $G$, or is a constant map
with image $a'$ equal to a fixed point of $G$.  

{\bf Case $H_\infty(W)=a'$ is constant}.  Then   
$G^r(c_\infty)=a'$ is a fixed point of $G$.  If $a'$ is attracting or
superattracting,  then the  basin of $a'$  contains $G^r(c_\infty)$. 
 By Proposition  \ref{prop:to_infty}, the basin of $a'$ is disjoint from 
$B_\infty$.   This is impossible, since then for $n$ large,
$F_{t_n}$ would have $F_{t_n}^r(c_{t_n})$ and $a_{t_n}$ in
distinct basins, a contradiction. Thus, $G^r(c_\infty)=a'$, a
repelling or neutral fixed point.  

{\bf Case $H_\infty|_W$ is univalent}. It follows easily that
$H_\infty(W)$ is contained  in  Fatou component $\Omega$ of $G$,
$G^r(c_\infty)\in \Omega$, and $G(\Omega)=\Omega$. We claim that
$\Omega$ must be parabolic; the condition $G(\Omega)=\Omega$ implies
that the multiplier is one. The same argument as that given in the
previous paragraph shows that $\Omega\ne B_\infty$ and  that $\Omega$
cannot be an attracting or superattracting basin. The map
$F|W$ has the property that $F(W)\subset W$ and that under iteration, every orbit leaves any
 compact
subset.  The same is true therefore for $G|_{H_\infty(W)}$. So $\Omega$ cannot
be a Siegel disk or Herman ring  either.  By the classification of Fatou
components, $\Omega$ is a parabolic basin.
\qed

{\noindent\bf Proof of Corollary} \ref{finite}. Now 
$F$ is hyperbolic, i.e. all critical points $c_i$ are in attracting
or superattracting basins. Assume, by taking a subsequence if
necessary, $ H_{t_n}(c_i)\to \xi_i$ as $n\to
\infty$. Clearly each $\xi_i$ is a critical point of $G$ and $G$ has
no other critical points. Recall that $B$ denotes the immediate
attracting basin containing the attractor
$a$ of the spun critical point.  Let $\widetilde{B}$ denote the full
basin of $a$.  

Suppose $c_i\not\in \widetilde{B}$.  Then $c_i\in \WWW_0$, where
$\WWW_0$ is as in Theorem \ref{thm:closure}. This theorem implies 
that $\xi_i=\JJJ(c_i)$ lies in an attracting or superattracing
basin of $G$.

Suppose $c_i\in \widetilde{B}$ and $c_i\ne c$. Again, Theorem
\ref{thm:closure} implies that $\xi_i=J(c_i)$ is attracted by
$a_\infty$.

Suppose $c_i=c$, the spun critical point.  Then $\xi_i=c_\infty$.
Theorem \ref{thm:coarse} implies that either
$G^r(c_\infty)$ is in a fixed parabolic basin $\Omega$,  or
$G^r(c_\infty)$ is a repelling or neutral fixed point. 

In the former  case, $c_\infty$ must be itself in $\Omega$ as well,
since $\Omega$ must contain a critical point, and by assumption
A3 the orbit of $c$ does not contain other critical points.
All other critical points converge to attractors, so $G$ has no other
parabolic basins and is therefore geometrically finite.

In the latter case, $G^r(c_\infty)$  cannot be neutral, since this
would require the existence of another critical point having infinite
orbit and not converging to an attracting cycle.
Hence $G^r(c_\infty)$ is a repelling fixed point, $G$ has no
parabolic basins, and $G$ is subhyperbolic, hence geometrically
finite.\qed

\section{Proof of Theorems 1.1, 1.2}
\label{secn:main}

Let $f$ be a critically generic hyperbolic rational map without
critical orbit relations, or a hyperbolic polynomial with connected
Julia set and without critical orbit relations .  Let
$a, B, \gamma, c, b$ be as in Theorem \ref{thm:lands_or_diverges}, and
suppose $c$ is visible after $r$ steps.  Let
$\sigma: [0, +\infty) \to \Rat_d/\Aut(\IP^1)$ be the spinning ray. 

In case that $f$ is a rational map, we make
the assumption that the spinning ray has at least one
limit point in $\Rat_d/\Aut(\IP^1)$.  In the polynomial case,
existence of a limit point follows since the connectedness locus is
bounded (Lemma \ref{lemma:conn_locus_bounded}).

If $f$ is a rational map, conjugate it and label its
critical points to produce $f^\times_0 \in \GRatmn$.  If $f$ is a
polynomial, use the orders of its critical points to determine a
partition $\DDD$ of $d-1$.  Conjugate $f$ and label its critical points
to produce $f^\times_0 \in \Poly^\times _d (\DDD)$.

Since $f$ is hyperbolic, and has no critical orbit relations, there
exists a decomposition of indices of critical points
\[ C=\{1, 2, ..., 2d-2\} = I \sqcup J \sqcup \{k\}\]
in the rational case and
\[ C=\{1, 2, ..., M\} = I \sqcup J \sqcup \{k\}\]
in the polynomial case, such that (i) $c_k=c$, the spun critical point,
and (ii) the conditions in \S 4 for the definition of the subspace $Z$
are satisfied.  By Corollary \ref{cor:one_dimensional}, the subspace
$X(f^\times_0)$ is a Riemann surface.  By
Proposition \ref{prop:lifts_to_X}, the spinning ray lifts to
$X(f^\times_0)$.

We now verify the assertion in Step 2 of the outline in \S 1.4.
Let $\Xi$ denote the set of limit points of the
spinning ray in $\Rat_d/\Aut(\IP^1)$, and let $\widetilde{\Xi}$ denote
the set of limit points of its lift in $X(f^\times_0)$.  Let
$\pi^{\times, *}: \GRatmn \to \Rat_d/\Aut(\IP^1)$ be the natural
projection which forgets the labelling of critical points and records
the conjugacy class.

Let $f^\times_t = \sigma^\times (t)$.  Suppose $\sigma(t_n) \to
\xi_\infty \in \Rat_d/\Aut(\IP^1)$.  We must show that after passing to
subsequences, $f^\times_{t_n} \to g^\times$ for some $g^\times \in
X(f^\times_0)$.  Since the possible labellings are finite in number, it
suffices to prove that the underlying maps $f_{t_n}$ converge to $g$.
The definition of the quotient topology on $\Rat_d/\Aut(\IP^1)$
implies that there exist $F_{t_n}, G \in \Rat_d$ such that $F_{t_n}
\to G$,  $\pi^{\times, *}(F_{t_n}) = \sigma(t_n)$,  and $\pi^{\times,
*}(G)=\xi_\infty$.  By Theorem \ref{thm:closure}, the critical points
of $G$ are distinct. We may  choose labelling of critical points so
that $F^\times_{t_n}
\to G^\times$ in $\GRatm$, since the set of possible labelling is
finite.  We may then write $f_{t_n}=M_n F_{t_n} M_n^{-1}$ where
$M_n(\{0,1,\infty\})$ is contained in the set of critical points of
$F_{t_n}$.  Since the critical points of $F_{t_n}$ converge to those of
$G$, which are distinct, after passing to a subsequence we have $M_n \to
M$.  Thus $f_{t_n} = M_n F_{t_n} M_n^{-1} \to MGM^{-1}$.
Putting $g=MGM^{-1}$, we have produced a limit point $g^\times$ of $f^\times_{t_n}$.
By Theorem \ref{thm:closure}, for $i\in I$,  the
attractors $a_i(g^\times)$ have the same multiplier,
the critical points $c_i(g^\times)$ have linearizing position
$-1$, and for $j\in J$, the critical points 
 $c_j(g^\times)$ have the same
linearizing coordinate as $c_j(f_0^\times)$. Thus $g^\times\in Z$. But
$X(f^\times_0)$ is  closed in $Z$. So $g^\times\in  X(f^\times_0)$,
and Step 2 is shown.  

By Corollary \ref{finite}, any limit point
$g^\times$ of $\sigma^\times (t)$ has either a $1$-parabolic
fixed point, or else $g^r(c_\infty)$ is a repelling fixed point.  That
is, $g \in X_{par}^{1, 1}$ or $g \in X_{mis}^{r, 1}$ in the notation
of \S 1.2.  Since
$f^\times_0$ is hyperbolic, 
 $X_{par}^{1,1}, X_{mis}^{r,1}$
are not all of $X(f^\times_0)$.     Since $X(f^\times_0)$ is
one-dimensional, it follows that $X_{par}^{1,1}, X_{mis}^{r,1}$
 are discrete subsets of $X(f^\times_0)$.  Since
the spinning ray is connected, if two distinct limit points exist, then a
continuum of limit points exists, which violates discreteness.
Hence the spinning ray $\sigma$ has a unique limit point in 
$ \Rat_d/\Aut(\IP^1)$.

This proves most of Theorems \ref{thm:lands_or_diverges}
and  \ref{thm:lands_for_polynomials}.
The proof of the conclusion that the limit is
independent of the representative $\gamma$
and the annulus $A$, requires a bit of technology from Teichm\"uller
theory and is given in \S 7.
\qed

\section{Limits of spinning, III}
\label{secn:fine}

In this section, we prove Theorems \ref{thm:fine_generic} and
\ref{thm:diverge}.  

We first recall some notation.  $A-\gamma$ is the union of
two annuli $A^\pm$, with core curve $\hat{\gamma}^\pm$ (see Figure
\ref{fig:torusplot.eps}). Denote by $\delta^{\pm}\subset
\Psi^{-1}(\hat{\delta}^\pm)$ (respectively,
$\gamma^{\pm}\subset \Psi^{-1}(\hat{\gamma}^\pm)$, $W\subset
\Psi^{-1}(A_c)$, $S^\pm\subset \Psi^{-1}(A^\pm)$
 ) the unique lifts  whose closures
contain the attractor $a$.  

Assume $c$ is visible.  Conclusion (1) follows
immediately from Corollary \ref{finite}.  

Taking subsequences if necessary, we may assume
$H_{t_n}|W\to H_\infty$ locally uniformly, $\overline{\gamma_{t_n}^\pm}\to \Gamma^\pm$
in the Hausdorff topology, $u_{t_n}^\pm\to u_\infty^\pm$, $a_{t_n}\to a_\infty$ and $c_{t_n}\to c_\infty$.
Note that $H_\infty(\gamma^\pm)\subset \Gamma^\pm$.
We will make use of the structure theorem for geometric limits of
invariant strips (Appendix \ref{app:geomlims}). 

In our setting $a_\infty\in \Gamma^\pm$ is an attracting vertex. It
follows by Theorem \ref{thm:structure} that $\Gamma^\pm$ has exactly
one edge $\zeta^\pm$ contained in $B_\infty$.
By Theorem \ref{thm:coarse} and  Corollary \ref{finite}, $c_\infty\in
H_\infty(W)\subset \Omega$, with $\Omega$ a fixed parabolic basin,
$H_\infty$ is univalent,
 $G$ has no other parabolic basins,
$\Omega$ contains no other critical points. 
Therefore $G|\Omega$ is conformally conjugate to the ``cauliflower''
$z^2+\frac14$. The set $H_\infty (W)$ represents an annulus in
$\Omega/G$ and contains the unique puncture (corresponding to
$c_\infty$).  Thus $H_\infty(\gamma^\pm)$ are two loops in
$\Omega$ ending at the parabolic point $v$, symmetric under an
anticonformal involution of $H_\infty(W)$.
By Theorem \ref{thm:structure}, every other edge of $\Gamma^\pm$ is
contained in $\Omega$, is
 therefore a loop based on $v$. So $\Gamma^\pm$ has only two vertices: $a_\infty$ and $v$. As 
the backward end of
 $\zeta^\pm$ is a non-attracting vertex of $\Gamma^\pm$, it must be $v$. So $v\in \partial
 B_\infty$.  In summary, $\Gamma^+\cup
\Gamma^-$ looks like Figure
\ref{fig:hlim.eps}, possibly with infinitely many loops.  This proves
conclusion (2).  

\xincludeps{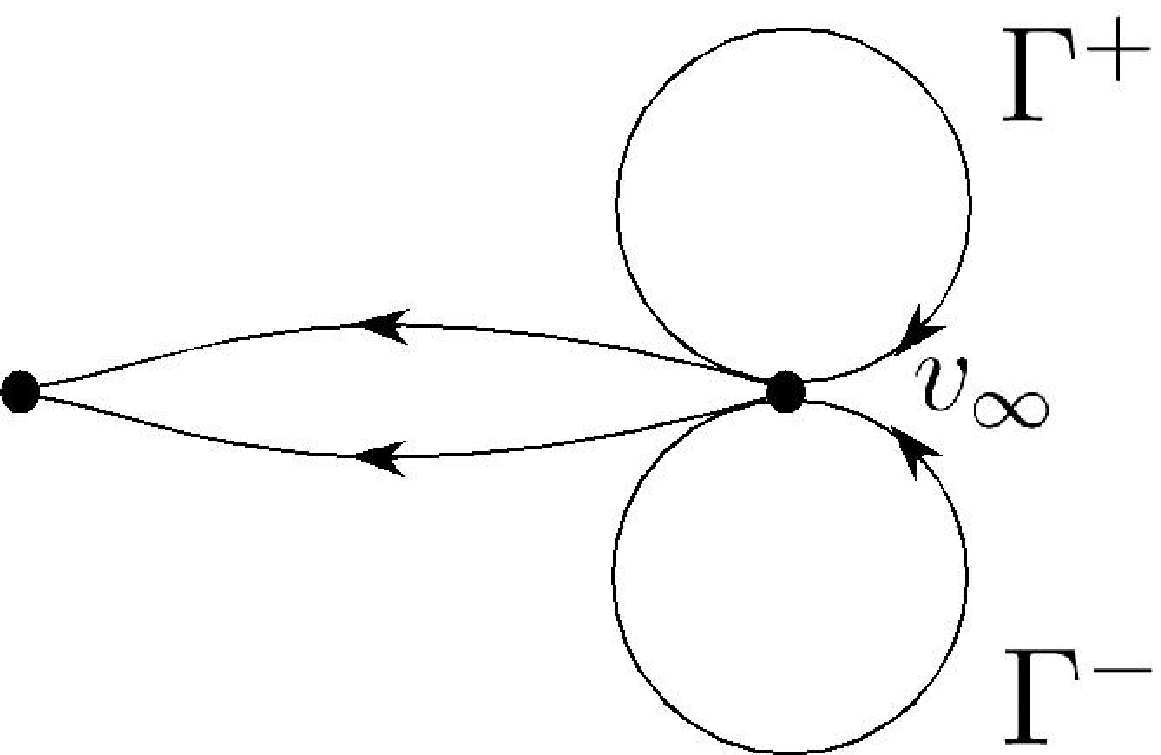}{2}{Structure of the Hausdorff limits $\Gamma^{\pm}$.}  

Conclusion (3) is more delicate.

Claim 1. The point $v$ split into two fixed points for nearby maps. This is due to the
fact that $v$ has multiplicity $1$. 

Denote by  $v_{t_n}^1, v_{t_n}^2$ the corresponding fixed points for $F_{t_n}$. As $F_{t_n}$ is
hyperbolic,  $v_{t_n}^1$,  $v_{t_n}^2$  are repelling, and thus distinct.

Claim 2.  $u_{t_n}^\pm\in \{v_{t_n}^1,v_{t_n}^2 \}$. This is because $u_\infty^\pm$ are 
non-attracting vertices of $\Gamma^\pm$,
must be equal to $v$ so
$u_{t_n}^\pm$ are close to $v$, and the only fixed points of $F_{t_n}$ close to $v$
are $v_{t_n}^1$ and $v_{t_n}^2$.

Choose now $P$ a repelling petal of $v$. $\partial P$ contains a subarc $I\subset B_\infty$ 
connecting $\zeta^+$ to $\zeta^-$, and $\partial P-(I\cup\{v\})$ 
has two components (call them 'sides'), each intersects one of $H_\infty(\gamma^+)$, $H_\infty(\gamma^-)$.

Claim 3. $\Gamma^+$ (resp. $\Gamma^-$) intersects $\partial P-\{v\}$ only 
on one side, e.g. the side of  $H_\infty(\gamma^+)$ 
(resp.  $H_\infty(\gamma^-)$).

Proof by contradiction. Take the outermost edges  $\iota_1,\iota_2$ of $\Gamma^+$ on each 
 side. Adjust $P$ so that each $\iota_i$ intersects $\partial P-\{v\}$ at only one point,
 transversally, and points outwards.
Enlarge $I$ into an arc $I'\subset \partial P$ so that it intersects only
 $\zeta^+$, $\iota_1$ and $\iota_2$ among the edges of $\Gamma^+$.
 Choose $U$ a small  neighborhood of 
$\overline{\zeta^+}$. Adjust $U$ and $P$ so that 
$\kappa=\partial (U-P)\cup I'$ is a $\Gamma^+$-transversal graph, intersecting 
$\Gamma^+$ only on $I'$, and at three points. By Theorem \ref{Transversal}, for large $n$, 
$\gamma_{t_n}^+\cap \kappa$ has essentially
the same oriented structure. However, 
$\gamma_{t_n}^+$ is an embedded arc without
self intersections. This is not possible by the Jordan curve theorem.

Claim 4. For some large $n$, $u_{t_n}^+\ne u_{t_n}^-$. 

Proof by contradiction. Denote by $\iota^\pm$ the outermost loop of
 $\Gamma^\pm-(\zeta^\pm\cup\{a_\infty\})$. Adjust $P$ and $I'$ as above.
Choose $U$ a small disk containing $\overline{\zeta^+}\cup\overline{\zeta^-}$. Define 
$\kappa$ as above.  Now $\Gamma^+\cup \Gamma^-$ intersects $\kappa$
at exactly four points, all pointing outwards. Again, for large $n$, 
$\gamma_{t_n}^+\cup \gamma_{t_n}^- \cap \kappa$ has essentially the
same oriented structure. 
Assume $u_{t_n}^+=u_{t_n}^-$. Then 
$\overline{\gamma_{t_n}^+}\cup \overline{\gamma_{t_n}^-}$ is a Jordan curve.
This is not possible.
(Note that the same idea proves also that $\zeta^+$ and $H_\infty(\gamma^+)$
are on the same side).

Claim 5.  $u^+\ne u^-$  (conclusion  (3)), as  $u_{t_n}^+\ne u_{t_n}^-$ for
some $n$, and $u_{t_n}^\pm=H_{t_n}(u^\pm)$. 

For the conclusion (4) we need the holomorphic index formula.
The only fixed points of $F_{t_n}$ near $v$ are $u_{t_n}^+,u_{t_n}^-$. If we
integrate around a small loop going once around $v$, then for $n$
sufficiently large, 
\[ \frac{1}{1-\lambda^+_{t_n}} + \frac{1}{1-\lambda^-_{t_n}}=
\frac{1}{2\pi i} \int \frac{dz}{z-F_{t_n}(z)} \to \frac{1}{2\pi i}\int 
\frac{dz}{z-G(z)}.\]
The first inequality now follows with $m$ equal to slightly less than
the real part of the index of $G$ at $v$. Since $|\lambda_{t_n}^\pm|>1$,
the second inequality is trivial.

Case $c$  not visible is much easier. We omit the details.

\qed 

Theorem \ref{thm:diverge} follows immediately from Theorem
\ref{thm:fine_generic}, since the  conclusion (3) (i.e. $u^+\ne u^-$)
of Theorem \ref{thm:fine_generic} is violated.

\section{Interpretation via Teichm\"uller theory}
\label{secn:teich}

In this section, we give an alternative construction of the complex
manifold $X(f^\times_0)$ presented in \S 4, and we prove the
uniqueness assertion of Theorems \ref{thm:lands_or_diverges} and
\ref{thm:lands_for_polynomials}. 
\gap

\noindent{\bf Moduli spaces.}  Fix 
\bi
\item a complex torus $T$
\item a nonempty set $\{z_1, ..., z_l\}$ of distinct points on
$T$
\item a nonzero primitive homology class $\alpha \in H_1(T,\Z)$.
\ib
We let $S=T-\{z_1, ..., z_l\}$, so that $\cl{S}=T$.  We recall that a
quasiconformal homeomorphism $h: S \to S'$ extends uniquely to a
quasiconformal homeomorphism $\cl{h}: T \to T'$.  

Recall that the {\em modular group} $\Mod(S)=QC(S)/QC_0(S)$, where
$QC(S)$ is the group of quasiconformal self-homeomorphisms of $S$ and
$QC_0(S)$ is the normal subgroup of those maps which are isotopic to the
identity through qc maps which leave the punctures fixed.   The group
$\Mod(S)$ (anti)-acts properly discontinuously by holomorphic
automorphisms on the Teichm\"uller space via  
\[ h.(\psi:S \to S') \equiv \psi \circ h^{-1}: S \to S'.\]

Define $P\Mod(S,\alpha)$ to be the subgroup of $\Mod(S)$
represented by those maps $h: S \to S$ for which $\cl{h}(z_j)=z_j,
j=1...l$ and for which $\cl{h}_*(\alpha)=\alpha$, where $\cl{h}_*:
H_1(T,\Z) \to H_1(T,\Z)$ is the induced map on homology.  It is a
subgroup of the pure modular group consisting of maps which fix each
puncture, but is not a normal subgroup.  

Let 
\[ \MMM(S,\alpha)=\Teich(S)/P\Mod(S,\alpha).\]

\begin{prop}
The space $\MMM(S,\alpha)$ is a complex manifold of dimension $l$.
\end{prop}

\pf It is enough to prove that the action is fixed-point free.   Let
$h \in P\Mod(S,\alpha)$, and let $(\psi: S \to S')$ represent a fixed
point for $h$.  Then there exists a conformal isomorphism $g: S' \to
S'$ such that $h$ is isotopic to $\psi^{-1} g \psi$.  Moreover,
$\cl{g}: T' \to T'$ fixes the primitive nonzero homology class
$\cl{\psi}_*(\alpha) \in H_1(T',\Z)$.  

By assumption $T'$ has at least one marked point $z_1$ which must
fixed by $\cl{g}$. Let $p:(\C,0) \to (T', z_1)$ denote a universal
cover, $\Lambda$ its deck group acting by translations, and let
$\tilde{g}: (\C, 0) \to (\C, 0)$ be a lift of $\cl{g}$ under $p$.   
Then $\tilde{g}(w) = \omega w$ for some $\omega \in \C^*$.  Since
$\cl{g}$ fixes a nonzero homology class, $\tilde{g}$ fixes a nonzero
element of $\Lambda$.  Hence $\omega = 1$, i.e. $\tilde{g}$ is the
identity.  Thus $g: S' \to S'$ is the identity and  
$h$ is isotopic to the identity.  
\qed

\gap
\noindent{\bf Puncture-forgetting maps.} As above, fix a torus $T$, a set 
$\{z_1, ..., z_l\} \subset T$, and let $S=T-\{z_1, ..., z_l\}$. Suppose
$l >1$ and write $l=n+k$ where $n, k > 0$.  Let $S^\sharp = T-\{z_1, ...,
z_n\}$.  

Suppose a point in $\Teich(S)$ is represented by $\psi: S \to S'$.  Let
 $S^{' \sharp} = T-\cl{\psi}(\{z_1, ..., z_n\})$ be the surface in which
the last $k$ punctures are filled in.  This induces a holomorphic
``puncture-forgetting'' map 
\[ \nu^\sharp: \Teich(S) \to \Teich(S^\sharp).\]

\begin{thm}
\label{thm:puncture_forgetting}
\be
\item The map $\nu^\sharp: \Teich(S) \to \Teich(S^\sharp)$ is a
holomorphic fibration; in particular it is a holomorphic submersion.

\item If $k=1$:
\be
\item the fiber $\FFF$ above a point $S^\sharp$ is a properly
embedded holomorphic disk which is naturally identified with the
universal  cover of $S^\sharp$, and  
\item there is a natural embedding 
\[ \pi_1(S^\sharp, z_{n+1}) \stackrel{\iota}{\hookrightarrow} \Mod(S)\]
such that the restriction of $\iota(\gamma)$ to $\FFF$ coincides with the
action of $\gamma$ on the universal cover of $S^\sharp$.  
\eb

\eb
\end{thm}
\pf See e.g. \cite{nag:book:teich}, \S 5.3 for (1) and \cite{kra:spins}
for (2). 
\qed

\noindent{\bf New construction of  $\mathbf{X(f^\times_0)}$}. We consider
only the case of rational maps; the polynomial case is entirely
analogous.   Assume a decomposition
$C=\{1, 2, 3, ..., 2d-2\}= I
\sqcup J \sqcup K $ and a function $\omega: J \to I$ are given as in \S 4,
and suppose $f^\times_0 \in Z$.  

Fix $i \in I$.   For $f^\times \in Z$, let $S_i(f^\times)$ denote the
quotient surface $\hat{B}_i/f$ corresponding to the $i$th attractor of
$f^\times$ (i.e. $S_i$ is the corresponding torus $T_i$ punctured at the
orbits of {\em all} critical points in the basin $\hat{B}_i$).  Let
$S^\sharp_i(f^\times)$ be the surface $S_i(f^\times) $ with the punctures
corresponding to free critical points $c_k, k \in K$ of $f^\times$ filled
in.  Let $\alpha_i(f^\times) \in H_1(T_i(f^\times), \Z)$ denote the
canonical homology class.   Abusing notation let us denote
by $c_i(f^\times), c_j(f^\times)$ ($j \in \omega^{-1}(i)$) both the
critical points and their images under projection to $T_i$.  

\begin{lemma}
The correspondence $f^\times \mapsto (S_i^\sharp(f^\times), 
\alpha_i(f^\times))$ determines
a well-defined holomorphic map $\varphi^\sharp_i: Z \to
\MMM(S^\sharp_i(f^\times_0), \alpha_i(f^\times_0))$.
\end{lemma}

\pf
Up to composition with an element of $P\Mod(S_i^\sharp(f^\times_0),
\alpha_i(f^\times_0))$, there is a unique isotopy class of quasiconformal
map $\psi: S_i^\sharp(f^\times_0) \to S_i^\sharp(f^\times)$ such that the
extension $\cl{\psi}: T_i(f^\times_0) \to T_i(f^\times)$ sends
$c_i(f^\times_0) $ to $c_i(f^\times)$,  sends $c_j(f^\times_0)$ to
$c_j(f^\times)$, and for which $\cl{\psi}_*(\alpha_i(f^\times_0)) =
\alpha_i(f^\times)$.  Hence $\varphi^\sharp_i$ is well-defined. 

To show that $\varphi^\sharp_i$ is holomorphic, use the fact that the
multiplier $\lambda_i(f^\times)$ and the locations $\phi_i(c_j)$ of the
critical points in the linearizing coordinates vary holomorphically in
$f$; see \S 4.  
\qed 

Denoting $\MMM(S_i^\sharp(f^\times_0), \alpha_i(f^\times_0))$ by
$\MMM^\sharp_i$, we have therefore a map
\[ \varphi^\sharp\equiv (\varphi^\sharp_i): Z \to \prod_i \MMM^\sharp_i.\]
 
\begin{lemma}
The map $\varphi^\sharp: Z \to \prod_i \MMM^\sharp_i$ is a holomorphic
submersion.
\end{lemma}

\pf Let $f^\times \in Z$. By (\cite{ctm:ds:qciii},
Thm. 6.2), the Teichm\"uller space of any rational map $f$ 
is naturally isomorphic to 
\[ \Teich(\Omega^{dis}/f) \times M_1(J,f) \times \Teich(\Omega^{fol},
f).\]
The second factor is a polydisk corresponding to invariant
line fields supported on the Julia set (conjecturally, this occurs only
in the case of Latt\`es examples).   The third is a polydisk
corresponding to deformations supported in Siegel disks, Herman rings,
and superattracting basins.  The first factor is in turn isomorphic to a
product of Teichm\"uller spaces of quotient surfaces.  Hence 
the Teichm\"uller spaces of quotient surfaces appear naturally as factors
in the complex manifold which is the Teichm\"uller space of $f$. 
Moreover, there is a canonically defined holomorphic map 
\[ \eta: \Teich(f) \to \Rat_d/\Aut(\IP^1)\]
obtained by straightening.  By taking the qc conjugacy to fix zero, one,
and infinity we get a lift $\eta^\times$ whose image lies in $Z$ by
construction. 

We have the following commutative diagram of pointed complex manifolds:

\renewcommand \arraystretch{2}
\[
\begin{array}{ccccc}
\prod_i\Teich(S_i(f^\times)) & 
                                        \stackrel{\iota}{\hookrightarrow} & 
                 \Teich(f) & \stackrel{\eta^\times}{\rightarrow} &
(Z,f^\times)\mbox{\hspace{.5in}}\\ (\nu^\sharp_i)=\nu^\sharp \downarrow &
\; &
\; &
\; &
\downarrow
\varphi^\sharp = (\varphi^\sharp_i)
\\
\prod_i\Teich(S^\sharp_i(f^\times)) & \; & \longrightarrow & \;
& \left( \prod_i \MMM^\sharp_{i}, (S_i^\sharp(f^\times))\right) 
\end{array}
\] 
Here, $\nu^\sharp_i: \Teich(S_i) \to \Teich(S^\sharp_i)$ is the map
induced by forgetting punctures corresponding to free critical points,
and $\iota$ is the inclusion map. 

The map on the bottom is a (universal) holomorphic covering map.  By Theorem
\ref{thm:puncture_forgetting},  each 
$\nu^\sharp_i: \Teich(S_i) \to \Teich(S^\sharp _i)$ is a holomorphic
submersion, therefore the product is a submersion as well.  Since the
diagram commutes, it follows that the derivative of $\varphi^\sharp$ is
surjective when evaluated at $f^\times$, and the proof is complete.  
\qed

Fix again $i \in I$ and let $n_i=\#\omega^{-1}(i)$.  Recall from \S 4
that 
\[ \mathbf{C}^i = (\C^*-\{(-1)\lambda_i^n,n\in \Z\})^{\omega^{-1}(i)} - \mbox{big diagonal}\]
where the big diagonal is the locus where two or more coordinates have the same $\lambda_i$-orbits.
The map $\varphi_i: Z \to
\MMM^\sharp_i$ can be written as a composition 
\[ Z \stackrel{(\lambda_i, \Phi_i)}{\longrightarrow} \Delta^* \times
\mathbf{C}^i \stackrel{\cl{\varphi_i^\sharp}}{\longrightarrow}
\MMM^\sharp_i.\]

The second map $\cl{\varphi^\sharp}_i$ is the one induced by sending the
pair 
$(\lambda, (w_1, ..., w_{n_i}))$ to the quotient torus
$\C^*/\genby{w \mapsto \lambda w}$ punctured at the images of
$-1$ and at the $w_j$'s.  Using this it is easy to see that the map 
$\cl{\varphi^\sharp}_i$ is in fact an infinite covering.  It follows that
the
$|K|$-dimensional analytic set  $X(f^\times_0)$ 
constructed in $\S 4$ coincides with the connected
component of $(\varphi^\sharp)^{-1}(\varphi^\sharp(f^\times_0))$
containing $f^\times_0$, which by the previous proposition is a complex
submanifold. 
\gap

\noindent{\bf When $\mathbf{|K|=1}$.} 
Suppose $|K|=1$ and the unique free critical point $c=c_k$ lies in
the basin of attraction of the $i$th attracting cycle of $f^\times_0$. 
Let $S_i=S_i(f^\times_0)$ denote the corresponding quotient surface,
regarded as the basepoint in
$\Teich(S_i)$.  Set $T_i=T_i(f^\times_0)$.  

The inclusion
\[ S_i \hookrightarrow S_i^\sharp = S_i \union \{c\} \]
induces the puncture-forgetting map 
\[ \nu^\sharp_i : \Teich(S_i) \to \Teich(S^\sharp_i).\]
By Theorem \ref{thm:puncture_forgetting}, the
fiber
$\FFF$ above the basepoint $S^\sharp_i$ is canonically identified with the
universal cover of $S^\sharp_i$ and is a properly embedded holomorphic
disk.  Spinning continuously about $\gamma$ defines a map $\hat{\sigma}:
\R \to \FFF$ such that $\hat{\sigma}(0)$ is the basepoint $S_i$. 
Identifying factors with their images under the canonical inclusion
(i.e. using basepoints corresponding to $f_0$ in the  factors other than
$i$) we have the following diagram of pointed manifolds:  

\[
\begin{array}{lcccl}
\; & (\R,0) &  \stackrel{\id}{\longrightarrow} & (\R,0) & \;\\
\; & \hat{\sigma} \downarrow & \; &  \downarrow \sigma^\times & \; \\
\Teich(S_i) \supset & (\FFF, S_i) & \stackrel{\eta^\times \circ
\iota}{\longrightarrow} 
              & (X(f^\times_0), f^\times_0) & \subset Z(f^\times_0) \\
\; & \nu^\sharp_i \downarrow & \; & \downarrow \varphi^\sharp_i& \; \\
\Teich(S^\sharp_i) \ni & \{S^\sharp_i\} &\longrightarrow 
&\{(S^\sharp_i, \alpha_i)\} & \in \MMM^\sharp_i  
\end{array}
\]

Let $t \in \Z$, let $i$ be the index of the attractor in which
spinning takes place, and let
$h: S_i \to S_i$ denote the quasiconformal homeomorphism as in the
definition of spinning (\S 1).  Then $h$ represents and element of
the pure modular group $P\Mod(S_i)$.  Since 
$\Mod(S_i)$ anti-acts on $\Teich(S_i)$ by {\em pre}composition on the
marking maps, if we set 
\[ \tau\equiv[h_{-1}]\]
then
\[ \tau.\hat{\sigma}(t)=\hat{\sigma}(t+1) \;\;\;\mbox{for all} \;\;\;
t \in \R.\]

The lemma below follows from Theorem \ref{thm:puncture_forgetting}.

\begin{lemma}
\label{lemma:hyperbolic}
The infinite cyclic group $\genby{\tau}$ acts on $\Teich(S_i)$ by
biholomorphic maps preserving the holomorphic disk $\FFF$.  With respect
to the Poincar\'e metric on $\FFF$, the map $\tau$ is a hyperbolic
translation whose length is the length of the unique simple closed
geodesic on $S^\sharp_i$ freely homotopic to $\gamma$.
\end{lemma}

To set up the next statement, let $H: \IP^1 \to \IP^1$ be a
quasiconformal conjugacy between $f^\times_0$ and another map
$F^\times_0 \in Z$ such that $c_m(F^\times_0) = H(c_m(f^\times_0))$,
$m=1, 2, ..., 2d-2$.  The pair $(F_0, H)$ represents an element of
$\Teich(f_0)$, by definition.

\begin{thm}
\label{thm:independence}
Suppose $(F_0,H) \in \FFF$.  Let $f^\times_n = (\eta^\times \circ
\iota )(\tau^n.(f_0,\id))$ and let $F^\times_n = (\eta^\times \circ
\iota )(\tau^n.(F_0,H))$.  
Then 
\[ \lim_{n \to \infty} F^\times_n  =
    \lim_{n \to \infty} f^\times_n\] if the latter limit exists in the
space $X(f^\times_0)$.   
\end{thm}

\pf Suppose the latter limit is $g^\times$.
Consider the hyperbolic surface $\Omega$ obtained by puncturing the
Riemann surface $X(f^\times_0)$ at $g^\times$ and equipping it with the
hyperbolic metric, denoted $d_\Omega$.   Let
$d_\FFF$ denote the hyperbolic metric on $\FFF$.  Since
$\tau$ acts by hyperbolic isometries, 
\[ d_\FFF(\tau^n.(f_0,\id), \tau^n.(F_0, H)) = d_\FFF((f_0,\id),
(F_0, H))=D,\]
which is independent of $n$.
By hypothesis,  the sequence $\{f^\times_n\}$ exits the cusp of
$\Omega$ corresponding to the puncture at $g^\times$. 
The map $(\eta^\times \circ \iota): \FFF \to \Omega$
is holomorphic, so it is distance non increasing with regard to $d_\FFF$
and $d_\Omega$.  Hence for all $n$, 
\[ d_\Omega(f^\times _n, F^\times _n) \leq D.\]
By choosing a local chart centered at $g^\times$ and comparing Euclidean
and hyperbolic metrics, we see that the sequence $\{F^\times_n\}$ must
exit the cusp at $g^\times$ as well.
\qed

\gap

\noindent{\bf Proof of uniqueness in Thms. 1.1, 1.2.}
Consider spinning, starting with the map $F^\times_0$, and using
the curve $H(\gamma)$ and the annulus $H(A)$ instead.  Then 
\[ \sigma_{\gamma, A}, \;\;\sigma_{H(\gamma), H(A)}: \R \to
\Rat_d/\Aut(\IP^1)\]
are two spinning paths, which yield lifts 
\[ \sigma^\times _{\gamma, A}, \;\;\sigma^\times _{H(\gamma),
H(A)}:
\R
\to X(f^\times_0)=X(F^\times_0).\]
Since 
\[ (\eta^\times 
\circ i)(\tau^n.(F_0,H))=\sigma^\times _{H(\gamma), H(A)}(n), \;\;\;n
\in \Z
\]
the previous theorem and Step 2 of the proof of Theorem 1.1 in \S 5 
yields 

\begin{cor}
The limits of spinning rays 
\[ \lim_{t \to +\infty} \sigma_{\gamma, A}(t) \;\;\mbox{ and } \;\;
\lim_{t \to +\infty} \sigma_{H(\gamma), H(A)}(t) \]
coincide if one (hence the other) exists.   In particular, the limit of
spinning depends only on the homotopy class of $\gamma$ in
$\pi_1(S^\sharp_i, c)$.  
\end{cor}

\section{Appendix:  geometric limits}
\label{app:geomlims}

In this section, we analyze the possible geometric limits of 
restrictions of rational maps to certain forward-invariant open
disks, called {\em invariant strips.}
\gap

\noindent{\bf Invariant strips.}  Let $L$ be a positive real number. 
The {\em standard strip $S$} is the domain 
\[ \{x+iy \; | \; |y| < L\}.\]
The {\em standard translation} $T: S \to S$ is given by $T(x+iy) = x+iy
- 1$, i.e. translation by one unit to the left.  The {\em standard
central line $\R$} is conformally distinguished as the unique 
geodesic (with respect to the hyperbolic metric on $S$) stabilized by
$T$.  In the following, we shall be concerned exclusively with the case when
the parameter $L$ is fixed.

\begin{defn}[Invariant strip]
Let $f: \IP^1 \to \IP^1$ be a rational map.  An 
{\em invariant strip of $f$} is an open,
simply-connected set $A
\subset \IP^1$ such that the restriction $f|A$ is holomorphically
conjugate to the standard translation on the standard strip.  The {\em
central line}
$\gamma$ of $A$ is the unique hyperbolic geodesic stabilized by $f|A$.  
\end{defn}

The following lemma follows from a normal families argument and the Snail
Lemma \cite{milnor:dynamics},  Lemma 16.2. 

\begin{lemma} [End points of an invariant strip] 
\label{lemma:endpoints} 
An invariant strip $A$ for a rational map $f$
has a unique {\em forward endpoint} $a$ given by $a=\lim_{n\to
+\infty} f^{\circ n}(z)$, where $z \in A$ is arbitrary.  $f(a)=a$,
and either $|f'(a)|<1$ or $f'(a)=1$. Similarly, it has a unique
fixed {\em backward endpoint}
$u$ where either $|f'(u)|>1$ or
$f'(u)=1$.  It is possible for the forward and backward ends to
coincide; this occurs if and only if $a=u$ is a 1-parabolic fixed
point of $f$. 
\end{lemma}

Assume, for the remainder of this section, that $f_n \to g$ uniformly,
$A_n \subset \IP^1$ are invariant strips for $f_n$ with central lines
$\gamma_n$, and $\cl{\gamma_n} \to \Gamma$ in the Hausdorff topology
of compact subsets. We now analyze the structure of $\Gamma$.

Denote by $a_n, u_n$ the  forward and backward ends of
$\gamma_n$, respectively.

\begin{thm}[Tameness of central lines]
\label{thm:tameness}
Suppose, taking subsequences if necessary, $a_n\to a_\infty$ and $u_n\to u_\infty$.
 Then, except when $\Gamma$ reduces to a single point,  there is a 
non-empty collection
$\{S_e\}_{e \in E}$ of disjoint invariant strips for $g$, indexed by a (finite or) countable
set $E$, with central lines $\gamma_e$, such that
\[ \Gamma-\Fix(g)=\bigcup_{e \in E} \gamma_e \quad \mbox{and}\quad \Gamma = \bigcup_{e \in E} \cl{\gamma_e}.\]
\end{thm}

\pf

Let $\HHH_n$ denote the set of all holomorphic conjugacies from $(T,S)$ to $(f_n,A_n)$ and put
$\HHH = \union \HHH_n$.  Then $\HHH$ is a normal family of univalent functions,
since the images of $S$ under elements of $\HHH_n$ must avoid three repelling
periodic points $x_n,y_n,z_n$ of $f_n$, with $x_n,y_n,z_n$ tending to $x,y,z$, three
distinct repelling periodic points of $g$.  We may assume one of these points, say $x_n$, is the point at
infinity.   Therefore, any sequence $(h_n)$ with  $h_n\in \HHH_n$ has a subsequence converging
uniformly on compact subsets.  Such a limit of univalent functions is either
univalent or constant. In either case, any limiting function $H$ satisfies
$H\circ T = g \circ H$. Moreover, suppose $h_{n_k}: (S,0) \to (A_{n_k}, y_k)$ is
any sequence of conjugacies which converges to a limiting conjugacy $h: (S,0) \to
(S_y, y)$.  Then $h$ is non-constant iff $g(y) \neq y$, and iff $ (A_{n_k}, y_k)\to (S_y,y)$
in the Carth\'eodory topology (see
e.g. \cite{ctm:renorm}, Ch. 5). 

\begin{lemma}
\label{lemma:same_or_disjoint}
Let $\{U_n\}$ be a sequence of open disks in $\C$ and let $z_n,z_n' \in U_n$.  
Suppose $(U_n, z_n) \to (U, a)$ and $(U_n, z'_n) \to (U', a')$ in the
Carath\'eodory topology.  Then either $U=U'$ or $U \intersect U' = \emptyset$.
\end{lemma}

Note that the basins $U_n$ are the same in both sequences; only the basepoints
differ. The proof is omitted.
\gap

\pf {\bf (of Thm., continued)}
Our goal is the construction of a sequence $(S_m, y_m)$ which
contains all possible
nonconstant limits of pointed strips $(A_n, y_n)$ where $y_n \in A_n$.

Choose a countable dense subset $\{y_m\}_{m=1} ^\infty$ of $\Gamma-\Fix(g)$.

Since $\overline{\gamma_n} \to \Gamma$ and $y_m \in \Gamma$, for each $m$, we may find a
sequence $\{y_{m,n}\}_{n=1}^\infty$ such that $y_{m,n} \in \gamma_n$
and $y_{m,n}
\to y_m$ as $n \to \infty$.

Consider now the following array:
\[
\begin{array}{cccc}
(A_1, y_{1,1}) & (A_2, y_{1,2}) & (A_3, y_{1,3}) & \ldots \\
(A_1, y_{2,1}) & (A_2, y_{2,2}) & (A_3, y_{2,3}) & \ldots \\
(A_1, y_{3,1}) & (A_2, y_{3,2}) & (A_3, y_{3,3}) & \ldots \\
\ldots & & &
\end{array}
\]
 From row 1, choose a subsequence $\{n(1,k)\}_{k=1}^\infty$ for which
the pointed
strips
$(A_{n(1,k)}, \; y_{1, n(1,k)})$ converge to a strip $(S_1, y_1)$.

 From row 2, using only those columns used in the subsequence just constructed,
choose a subsequence $\{n(2,k)\}_{k=1}^\infty$ for which the pointed strips
$(A_{n(2,k)}, \; y_{2, n(2,k)})$ converge to a strip $(S_2, y_2)$.

Continue this process inductively, each time using only column
indices that have
already been used.  We obtain, for each row $m$, a sequence
$\{n(m,k)\}_{k=1}^\infty$ of column indices such that
\[ (A_{n(m,k)}, \; y_{m, n(m,k)}) \to (S_m, y_m).\]

This construction has the following property.  Suppose $m <  m'$,
i.e. that row $m$
is above row $m'$. Then
\begin{equation}
\label{eqn:m}
(A_{n(m,k)}, \; y_{m, n(m,k)}) \to (S_m, y_m)
\end{equation}
while also
\begin{equation}
\label{eqn:mp}
(A_{n(m',k)}, \; y_{m', n(m',k)}) \to (S_{m'}, y_{m'}).
\end{equation}
However, by construction the sequence $\{n(m',k)\}_{k=1}^\infty$ of
column indices appearing in row $m'$ is in fact a subsequence of
$\{n(m,k)\}_{k=1}^\infty$, the  indices for  row $m$, so we may
substitute $m=m'$
in the expression $n(m,k)$ appearing in equation (\ref{eqn:m}) to deduce that
\begin{equation}
\label{eqn:msub}
(A_{n(m',k)}, \; y_{m, n(m',k)}) \to (S_m, y_m).
\end{equation}
By Lemma \ref{lemma:same_or_disjoint} applied to the two converging
pointed strips
in (\ref{eqn:mp}) and (\ref{eqn:msub}), we conclude that either $S_m =
S_{m'}$ or else $S_m$ and $S_{m'}$ are disjoint.

Now suppose $y \in \Gamma - \Fix(g)$.  We will show $y$ lies on a central line
$\gamma_m$ of $S_m$ for some
$m$. Choose a subsequence $\{y_{m_i}\}_{i=1}^\infty$ with $y_{m_i} \to y$ as $i \to
\infty$.  Then from the array shown above we may extract a diagonal
sequence of pointed basins
\begin{equation}
\label{eqn:diag}
(A_{n(m_i, k_i)}, \; y_{m_i, n(m_i, k_i)}) \to (S_y, y)
\end{equation}
where $S_y$ is  an invariant  strip containing $y$.

We will now  show that $S_y$ is one of the limiting invariant  strips already
constructed.  Fix a row, say $m$.  By construction, the  sequence  of
column entries
occurring in  (\ref{eqn:diag}) is eventually a subsequence of the column entries
chosen for row $m$, i.e. $n(m_i, k_i) \in \{n(m,k)\}_{k=1}^\infty$ for all $i$
sufficiently large.  Hence in addition  to (\ref{eqn:diag}) we may obtain, by
substituting $n(m,k)=n(m_i,  k_i)$ in (\ref{eqn:m}) that
\begin{equation}
\label{eqn:diagsubs}
  (A_{n(m_i, k_i)}, \; y_{m, n(m_i, k_i)}) \to (S_m, y_m).
\end{equation}
Comparing (\ref{eqn:diag}) and (\ref{eqn:diagsubs}) we see that Lemma
\ref{lemma:same_or_disjoint} applies again.  So, for any fixed $m$, either $S_m = S_y$ or else they
are disjoint.

However, $S_y$ is open and $\{y_m\}$ is dense in $\Gamma-\Fix(g)$, so $S_y$ cannot be
disjoint from all
of the  strips $S_m$.  Hence $y \in  S_m$ for some $m$.  Since $y$ is on the
central line $\gamma_y$ for $S_y=S_m$ and central lines are unique,
we have $y \in
\gamma_m$. We now enumerate the distinct central lines which arise as
$\{\gamma_e\}_{e \in E}$. Clearly 
$$\Gamma-\Fix(g) = \bigcup_{e\in E}\gamma_e.$$

In case $\Gamma$ reduces to a single fixed point $E$ is empty. Otherwise $\Gamma$ is compact, 
connected, and contains
infinitely many points. Note that  $\Fix(g)$ is finite.

We now show $\Gamma=  \bigcup_{e\in E}\overline{\gamma_e}$ by contradiction:

Assume $y\in \Gamma\cap \Fix(g)$ and $y\notin\overline{ \gamma_e}$ for any $e$. Choose a small 
closed disk $D$ centered at $y$
such that $D$ does not contain other fixed points of $g$, and $\Gamma-D\ne \emptyset$.

There are only finitely many edges $\gamma_i$ intersecting $\partial D$ (as they
belong to disjoint strips). If $y\notin\overline{\gamma_i}$, then
the ends of $\gamma_i$, as fixed points of $g$, are outside $D$. Using
the fact that $\Gamma$ is connected, one gets easily a contradiction.
\qed

The previous result does not exploit the constraints imposed by plane
topology.  A $C^1$-smooth arc $\kappa$ is called a $\Gamma$-{\em
transversal} if its ends do not meet $\Gamma$, it does not meet the
vertices $V$  and it intersects any edge transversally (if not
empty). One defines a $\Gamma$-{\em transversal} graph (with finitely
many arcs) similarly.

\begin{thm}[Stability of transversal]\label{Transversal}
Let $\kappa$ be a $\Gamma$-transversal graph. Then taking a subsequence if necessary, 
\be
\item $\Gamma\cap \kappa$ is finite.
\item There is $n'$, for any $n\ge n'$, $\#\overline{\gamma_n}\cap \kappa=\#\gamma_n\cap
 \kappa=\#\Gamma\cap \kappa$,
$\overline{\gamma_n}\cap \kappa \to \Gamma\cap \kappa$ in the Hausdorff topology.
Furthermore, $\gamma_n$ is transversal to $\kappa$, in the same orientation
as $\Gamma\cap\kappa$.
\eb
\end{thm}
\pf
Each $\gamma_e$ is a real-analytic arc. By transversality it meets $\kappa$ at at most finitely many points.
Only finitely many edges of $\Gamma$ can meet $\kappa$. So $\Gamma\cap \kappa=\bigcup\{y_i\}$ is finite.

For each $i$ 
choose  $y_{i,n}\in \gamma_n$ converging to $y_i$. Let $h_{i,n}: S\to A_n$ be a conjugacy
 from the standard translation
$T$ to $f_n$, mapping $0$ to $y_{i,n}$. Taking subsequences if necessary, we may assume, for every $i$,
$h_{i,n}|_S\to l_i|_S$ locally uniformly. As $l_i(0)=y_i$ and $y_i$ is not fixed by $g$, $l_i$ is univalent.
Therefore $l_i(\R)$ is an edge of $\Gamma$.

Cover now $\Gamma\cap \kappa$ by finitely small disks $\bigcup D_i$, with boundary
 transversal to both $\Gamma$ and $\kappa$.
We may choose $D_i$ small enough so that  $D_i\subset l_i(\Delta)$,
with $\Delta$ some fixed small closed neighborhood of $0$.
For $n$ large, $\overline{\gamma_n}\cap \kappa \subset \bigcup D_i$ by Hausdorff convergence. 
Note that, as analytic functions, the derivatives $h_{i,n}'$ converges to $l_i'$ locally uniformly as well,
in particular uniformly on $\Delta$. By transversality with $\kappa$, $(h_{i,n}(\R\cap \Delta))\cap \kappa$ 
is a single point, is contained in $D_i$, and the intersection is transversal,
 in the same direction as $l_i(\R)$.
Now $\gamma_n= h_{i,n}(\R)$ can not cross $D_i$ again, as $D_i\subset h_{i,n}(\Delta)$ and $h_{i,n}$ is
univalent.
\qed

\begin{thm}[Structure theorem for limits of strips]
\label{thm:structure} 
Under the hypotheses of Theorem \ref{thm:tameness}, 
the set $\Gamma$ admits the structure of a connected directed
planar graph
$\Gamma = (V,E)$ (possibly with
edges joining a vertex to itself, and possibly  with infinite
valence), where the
edges point in  the direction of the dynamics, subject to the following
restrictions:
\be
\item $V = \Gamma \intersect \Fix(g)$ is nonempty and finite, and contains $a_\infty,u_\infty$.

\item $\Gamma-V= \bigcup_{e\in E}\gamma_e$. Each edge $\gamma_e$ 
is the central line  of an invariant strip $S_e$ (recall that
these strips are disjoint), so
\item the edges are isolated
\item every vertex is the end of at least one edge
%\item between any two distinct vertices, there are at most finitely many edges

\item  $a_\infty$ is either an attracting or a parabolic fixed point of $g$. In the former
 case, $u_\infty\ne a_\infty$ and
there is a unique edge pointing to $a_\infty$ (and no edges pointing out from $a_\infty$). 
(there is a symmetric statement for $u_\infty$).
\item if the forward end of an edge $\gamma_e$ is not $a_\infty$, then $S_e,\gamma_e$ are
 contained in a fixed parabolic
basin.

\eb
\end{thm}

\pf 
Assume that $w\ne a_\infty$ is the forward end of a  $\gamma_e$. By
the Snail Lemma $w$ is either attracting or $1$-parabolic. It cannot 
be attracting as otherwise  by stability some point of $\gamma_n$ 
would be in a distinct attracting basin than that of $a_n$. So $w$ is
$1$-parabolic.

Assume $a_\infty$ is attracting and there are more than one edges  ending at $a_\infty$.
Let $D$ be a small disk centered at $a_\infty$, contained in the attracting basin, and whose
boundary is a $\Gamma$-transversal. 
By Theorem \ref{Transversal}, for large $n$,  $\gamma_n$ intersects transversally $\partial D$ with
the same number of points, all pointing inside $D$. But $\gamma_n$ is the central line
of a single invariant strip. This is not possible by the stability of the attractors.
\qed

\section{Appendix:  analytical lemmas}
\label{app:analytical}

\begin{lemma}
\label{lemma:koenigs_cover}
Let $B$ the full basin of attraction of a fixed point $a$ of
multiplier $\lambda$ with $0 < |\lambda| < 1$ for a rational map $f$, and
let $\psi: (B,a) \to (\C,0)$ be a  linearizing function for $B$.
Then:
\be
\item the set of critical points $\Crit(\psi)$ of $\psi$ is
$ \bigcup_{n \geq 0} f^{-n}(\Crit(f) \intersect B)$
and is invariant under $f^{-1}$.
\item the set of critical values of $\psi$ is
$\psi(\Crit(\psi)) = \bigcup_{n \geq 0} \lambda^{-n}\psi(\Crit(f)
\intersect B)$
and is invariant under multiplication by $\lambda^{-1}$.
\item if $B' = B - \Crit(\psi)$, then
$ \psi: B' \to \psi(B')$
is a covering map onto its image.
\eb
\end{lemma}
\pf
These properties follow immediately from the existence of a local
biholomorphic conjugacy on a neighborhood of the attractor $a$, and the
extension of this conjugacy to $B$ by setting $\psi(z) =
\lambda^{-1}\psi(f(z))$.
\qed

Recall that a {\em holomorphic family of rational maps over $X$} is a
holomorphic
map
\[ f: X \times \IP^1 \to \IP^1\]
where $X$ is a connected complex manifold.  

\begin{lemma}[Holomorphic dependence of K\"onigs functions]
\label{lemma:holo_koenigs}
Let
\[ f(x, z) = f_x(z)=\lambda(x)z + a_2(x)z^2 + ...\]
be a holomorphic family of rational maps over $X$ for which
$0<|\lambda(x)| < 1$,  and let $g \in X$.  Then there exists a neighborhood
$\UUU$ of $g$ and a disk
$D \subset \C$ containing the origin such that
\be
\item for all $f \in \UUU$, we have $f(\cl{D}) \subset D$, $f|D$ is univalent,
and
$D \subset B_f(0)$, the immediate basin of attraction of the origin for $f$
\item the map $\psi: \UUU \times D \to \C$ given by $\psi(x,z) =
\psi_{f_x}(z) = z + b_2(x)z^2 + ...$ (where $\psi_{f_x}$ is the
unique linearizing function for
$f_x$ normalized to have derivative one at the origin),  is holomorphic.
\eb
Furthermore, given any compact subset $K$ of $B_g(0)$, there is a neighborhood
$\UUU$ of $g$ in $X$, a neighborhood $\NNN$ of $K$ in $B_g(0)$ containing the
origin, and an extension of $\psi$ to a holomorphic map $\psi: \UUU \times \NNN
\to \C$.
\end{lemma}

The proof is omitted. 
 There is also a similar statement in the
superattracting case, assuming that all $f(x,z)$ have the same order
at $0$.

\begin{lemma}[Limits of models]
\label{lemma:limits_of_models}
Suppose $f, f_n, g$ are rational maps with $f_n \to g$ uniformly on $\IP^1$.
Suppose $W$ is a connected open hyperbolic subset of $\IP^1$, with
 $f^k(W) \subset W$, and suppose $j_n: W \to W_n$ are
holomorphic embeddings such that $j_n \circ f^k = f^k_n \circ j_n$.
Then, after passing to subsequences, the maps $j_n$ converge  uniformly on
compact subsets of $W$ to a holomorphic map $j$ such that $j \circ f^k = g^k \circ
j$.  Moreover, either $j$ is univalent, or $j(W)$ is a single fixed point of $g^k$. 
\end{lemma}
\pf
 Let $x,y,z$ be three distinct repelling periodic points of $g$ of periods
$p_x, p_y, p_z$.  By the stability of repelling periodic points,  we may find
repelling periodic points $x_n, y_n, z_n$ of $f_n$ of periods $p_x, p_y, p_z$,
respectively,  such that $x_n, y_n, z_n \to x, y, z$ as $n \to \infty$.   Since
$j_n(W) \subset \IP^1-\{x_n, y_n, z_n\}$, we have that $\{j_n\}$ is a normal
family.  Thus after passing to subsequences we have
$j_n \to j$ uniformly on compact subsets. Moreover $j \circ f^k = g^k \circ
j$. \qed

\begin{lemma}[Partial persistence of the dynamics]
\label{lemma:persistence}
Suppose $f, f_n, g$ are rational maps with $f_n \to g$ uniformly on $\IP^1$,
and with $f_n = h_n f h_n^{-1}$, where $h_n: \IP^1 \to \IP^1$ are  homeomorphisms.
Suppose $W$ is a periodic Fatou component of $f$ that is not a Siegel disk, and suppose that
%if $W$ is a Siegel disk, assume that the rotation number is Bryuno.  
$h_n$ are holomorphic on
the grand orbit $\widetilde{W}$ of $W$.

Then after passing to subsequences, the maps $j_n = h_n|\widetilde{W}$
converge uniformly on compact subsets to a univalent $j: \widetilde{W} \to
\IP^1$ satisfying $j \circ f = g \circ j$. 
If $W$ is attracting or superattracting, then $j(W)$ is a Fatou
component of $g$.  In general, $j(W)$ is containd in the Fatou set of
$g$.
\end{lemma}

\pf Denote by $\LLL$ the finite set of critical points and neutral periodic points. For any $w\in \LLL$,
assume (by taking subsequences if necessary) $h_n(w)\to w_\infty$.
If $w$ is a critical point, so is $w_\infty$, possibly with higher order.
If $w$ is eventually periodic, so is $w_\infty$.
If $w$ is $p$-periodic point of multiplier $\lambda\in S^1-\{1\}$, so is $w_\infty$.
If $w$ is of multiplier $1$, so is $w_\infty$, but with possibly a smaller period. 

On a given component, $j$ is either univalent or constant.  By the functional
equation, it suffices to show that $j$ is
univalent on periodic Fatou components.   By considering iterates, it
is enough to
prove that $j|W$ is nonconstant, where $W$ is a forward-invariant
component of $f$.

Case (a): $W$ contains a superattracting fixed point $0$ of order $k$. Then $f_n$ 
has a superattracting fixed point at $h_n(0)$ of the same order. Assuming $h_n(0)\to 0$. Then
$0$ is a superattracting fixed point of $g$, possibly with a higher order.

Assume by contradiction that $j|W$ is constant. Then $B_g(0)$ has no other critical 
points then $0$, for if $c_g\in B_g(0)$, and $c_g\ne 0$, then there are $c_n$ critical point
of $f_n$ contained in $h_n(W)$ and away from $h_n(0)$. There
is therefore $c\in W$ such that $j_n(c)=c_n$ for infinitely many $n$. This contradicts
 $j|W\equiv 0$.

Therefore $B_g(0)$ is simply connected.

Take a large disk $K\subset B_g(0)$ such that $g(K)\subset interior(K)$ and $K$ contains
all critical points in $B_g(0)$. This is stable for $f_n$. Therefore $h_n(W)$ is simply connected.
So is $W$. This then in turn contradicts the fact that  $j_n$ must preserve the modulus of
$W-L$ for any $L$ compact in $W$. 

Case (b):  $W$ contains an attracting but not superattracting fixed point $a$.  Then
$a_n=j_n(a)$ is an attractor for $f_n$ with the same multiplier, and after
passing to a subsequence we have
$a_n \to a$, an attractor for $g$ of the same multiplier. But $W$ contains a critical
point $c'$. So if $j$ were constant $j(a)=j(c')$ must be a superattracting fixed point.
This is not possible.

Case (c):  $W$ is a parabolic basin.  Then $W$
contains a critical point $c'$.  If $j$ is constant, then after passing to
subsequences we have as $n \to \infty$ that $c'_n = j_n(c') \to c'_g$ , a fixed
critical point of $g$.  But then for $n$ sufficiently large, the critical point
$c'_n$ converges under iteration of $f_n$ to an attracting fixed point $a_n$,
which is not the case.  

%Case (c):  $W$ is a Siegel disk.  Let $a_n = j_n(a)$ denote the
%centers.  If $j$ is constant,  then after passing to subsequences we have $a_n
%\to a \in \Fix(g)$, where the multipliers are constant.   By Bryuno's Theorem
%(\cite{milnor:dynamics}, Thm. 11.10),  $g$ has a Siegel disk centered at $a$ as well.  The proof
% concludes as in case (a).

Case (d):  $W$ is a
Herman ring.    Let $K$ be a leaf of the foliation of the Herman
ring.   If $j$ is
constant, then the spherical diameter of $j_n(K)$ tends to zero.
Hence the diameter
of a component $D_n = h_n(D)$ of the complement of $\gamma$ tends to
zero as well.
But $D \intersect J \neq \emptyset$ and this yields a contradiction
to Lemma 3.1 of
\cite{tan:pinching}. The same lemma yields also that $j(W)$ is in fact contained in a Herman ring
(but not a Siegel disk).
\qed

{\noindent\bf Remark.} A. Ch\'eritat (\cite{Che}, Part II) proved in a more general setting that
 if $W$ is a Siegel disk (resp. Herman ring) 
of Brjuno rotation number, then $j|_{\widetilde{W}}$ is univalent and $j(W)$ is again a Siegel disk (resp. Herman ring).

\newpage
\bibliographystyle{math}

\noindent \textsc{Kevin M. Pilgrim\\  
Dept. of Mathematics\\
Rawles Hall\\
Indiana University\\
Bloomington, IN 47405-7106, U.S.A. }\\
pilgrim@indiana.edu
\gap\gap

\noindent \textsc{Tan Lei\\
Unit\'e CNRS-UPRESA 8088\\
D\'epartement de
Math\'ematiques\\
Universit\'e de Cergy-Pontoise\\
2 Avenue Adolphe Chauvin\\
95302 Cergy-Pontoise cedex, France} \\
tanlei@math.u-cergy.fr

\end{document}